\documentclass[12pt,draftcls,onecolumn]{IEEEtran}

\usepackage{amssymb,amsmath,amsfonts,amsthm}
\usepackage{algorithm,algorithmic}
\usepackage{boxedminipage}
\usepackage{cite}
\usepackage[dvips]{graphicx}
\usepackage{epsfig}
\usepackage{float}
\graphicspath{{figures/}}
\usepackage{graphicx}
\usepackage{graphics}
\usepackage{multirow}
\usepackage{threeparttable} 
\usepackage[usenames,dvipsnames]{color}
\usepackage{url,bm,xspace,dsfont}

\usepackage{verbatim}

\def\diag{\mathop{\mathrm{diag}}}  

\newtheorem{theorem}{Theorem}
\newtheorem{lemma}{Lemma}
\newtheorem{definition}{Definition}
\newtheorem{corollary}{Corollary}

\newtheorem{remark}{Remark}

\newtheorem{assumptions}{Assumptions}
\newtheorem{problem}{Problem}

\newcommand{\eps}{\epsilon}

\newcommand{\veps}{\varepsilon}
\newcommand{\la}{\langle}
\newcommand{\ra}{\rangle}

\newcommand{\sr}{\stackrel}

\newcommand{\rar}{\rightarrow}

\newcommand{\tri}{\sr{\triangle}{=}}

\newcommand{\be}{\begin{equation}}
\newcommand{\ee}{\end{equation}}
\newcommand{\bea}{\begin{eqnarray}}
\newcommand{\eea}{\end{eqnarray}}
\newcommand{\bes}{\begin{eqnarray*}}
\newcommand{\ees}{\end{eqnarray*}}

\newcommand{\bi}{\begin{itemize}}
\newcommand{\ei}{\end{itemize}}
\newcommand{\ben}{\begin{enumerate}}
\newcommand{\een}{\end{enumerate}}

\newcommand{\bp}{\begin{problem}}
\newcommand{\ep}{\end{problem}}
\newcommand{\hso}{\hspace{.1in}}
\newcommand{\hst}{\hspace{.2in}}

\newcommand{\noi}{\noindent}

\newcommand{\bc}{\begin{center}}
\newcommand{\ec}{\end{center}}

\floatstyle{ruled}
\newfloat{algorithm}{H}{algo}
\floatname{algorithm}{\footnotesize Algorithm}

\begin{document}
%
\title{\bf Team Games Optimality Conditions  of Distributed Stochastic Differential Decision Systems with Decentralized Noisy Information Structures}


\author{ Charalambos D. Charalambous\thanks{C.D. Charalambous is with the Department of Electrical and Computer Engineering, University of Cyprus, Nicosia 1678 (E-mail: chadcha@ucy.ac.cy).}  \: and Nasir U. Ahmed\thanks{N.U Ahmed  is with the School of  Engineering and Computer Science, and Department of Mathematics,
University of Ottawa, Ontario, Canada, K1N 6N5 (E-mail:  ahmed@site.uottawa.ca).}
}

\maketitle

\begin{abstract}

We consider a team game reward,  and we derive a stochastic Pontryagin's maximum principle for distributed stochastic differential systems  with decentralized noisy information structures.  Our methodology utilizes the semi martingale representation  theorem, variational methods, and backward stochastic differential equations.  Furthermore, we derive necessary and sufficient optimality conditions that characterize team and person-by-person  optimality of  decentralized strategies.

Finally, we apply the stochastic maximum principle to several examples from the application areas of communications, filtering and control.
\end{abstract}

\vspace*{1.5cm}

  \vskip6pt\noindent {\bf Index Terms.}  Team Games, Distributed Systems, Decentralized, Optimality Conditions, Maximum Principle.

  \section{Introduction}
\label{introduction}
We derive necessary and sufficient team game optimality conditions for distributed stochastic differential systems with decentralized noisy information structures. For noiseless information structures, analogous optimality conditions are derived recently in \cite{charalambous-ahmedFIS_Parti2012} utilizing the representation of Hilbert space semi martingales and the stochastic Pontryagin's maximum principle of partially observed stochastic differential systems developed in \cite{ahmed-charalambous2012a}.  However, the   results obtained in \cite{charalambous-ahmedFIS_Parti2012} for decentralized noiseless information structures are not necessarily applicable to decentralized noisy information structures.
 In fact, there are certain technicalities that must be addressed when dealing with noisy information structures, which are inherited from the centralized fully observable versus partially observable stochastic optimal control  \cite{fleming-rischel1975,elliott1977,bismut1978,elliott1982,elliott-kohlmann1994,peng1990,yong-zhou1999,bensoussan1992a,charalambous-hibey1996,ahmed-charalambous2007}. The main  underlying assumption for centralized information structures, is that the acquisition of the information is centralized or the information acquired at different  locations  is communicated to each decision maker or control.\\

When the system model consist of multiple decision makers,  and  the acquisition of information and its processing is decentralized or shared among several locations, then the different  decision makers  actions are based on different information \cite{varaiya-walrand1978}.  We call  the information available for such  decisions,  "decentralized information structures or patterns"  \cite{witsenhausen1968,witsenhausen1971}.  When the system model is dynamic, consisting of an interconnection of at least two subsystems, and  the decisions are based on decentralized information structures, we call the overall system  a "distributed system with decentralized information structures".\\

Over the years several specific forms of decentralized information structures   are analyzed mostly in discrete-time 
\cite{witsenhausen1968,witsenhausen1971,ho-chu1972,kurtaran-sivan1973,sandell-athans1974,kurtaran1975,varaiya-walrand1978,ho1980,bagghi-basar1980,krainak-speyer-marcus1982a,krainak-speyer-marcus1982b,bansal-basar1987,waal-vanschuppen2000}, and more recently \cite{bamieh-voulgaris2005,aicardi-davoli-minciardi1987,nayyar-mahajan-teneketzis2011,vanschuppen2011,lessard-lall2011,mahajan-martins-rotkowitz-yuksel2012,gattami-bernhardsson-rantzer2012,mishra-langbort-dullerud2012,farokhi-johansson2012}.  
 However,   at this stage the only systematic framework  addressing  optimality conditions for distributed systems with decentralized information structures is the one reported in \cite{charalambous-ahmedFIS_Parti2012} for decentralized noiseless information structures. \\

 In this paper, we consider a team game reward  \cite{marschak1955,radner1962,marschak-radner1972,krainak-speyer-marcus1982a,waal-vanschuppen2000}, and we derive necessary and sufficient optimality conditions for distributed stochastic differential systems  with decentralized noisy information structures.        Our methodology utilizes the semi martingale representation  theorem, variational methods, and generalizes the concepts utilized in \cite{ahmed-charalambous2012a,charalambous-ahmedFIS_Parti2012} to derive  optimality conditions for   nonlinear stochastic distributed systems with decentralized  noiseless information structures.  From the practical point of view, the results of this part give optimality conditions in terms of forward and backward stochastic differential equations, and a Hamiltonian, called "Hamiltonian System of Equations", which we use to compute the optimal decentralized decision strategies of several examples from the application areas of communications and control.\\

The specific objectives of this paper  are the following.

\textbf{(a)} Derive team games Pontryagin's stochastic minimum principle (necessary conditions of optimality) for distributed stochastic  systems with decentralized noisy information structures;

\textbf{(b)} Introduce assumptions so that the necessary conditions of optimality in {\bf (a)} are also sufficient, and relate the optimality conditions to person-by-person optimality conditions;

\textbf{(c)} Apply the stochastic minimum principle to several examples from the application areas of communication and control.

The rest of the paper is organized as follows. In Section~\ref{formulation}  we  we formulate the distributed stochastic system with decentralized information structures.  In Section~\ref{sexistence}, we introduce the variational equation and discuss its application in decentralized filtering and control. Section~\ref{smp}  is devoted to the development of stochastic optimality conditions for team games with decentralized information structures,  consisting of necessary and sufficient conditions of optimality. In Section~\ref{classes} we apply the minimum principle to various examples.    The  paper is concluded with some comments on possible extensions of our results.

 \section{Distributed Stochastic Differential Team Games}
 \label{formulation}

 In this section we introduce  the mathematical formulation of distributed stochastic  differential systems, the noisy information structures available to the decision makers,  and the definitions of collaborative decisions via team game optimality and person-by-person optimality. Although, the stochastic differential systems are driven by the Decision Makers  (DMs) actions, our analysis includes unforced stochastic differential systems modeling distributed estimation.  Therefeore, the term "decision maker" is used for distributed control as well as distributed estimation.\\

The formulation presupposes a fixed probability space with filtration,   $\Big(\Omega,{\mathbb F},  \{ {\mathbb F}_{0,t}:   t \in [0, T]\}, {\mathbb P}\Big)$  satisfying the usual conditions, that is,  $(\Omega,{\mathbb F}, {\mathbb P})$ is complete, ${\mathbb F}_{0,0}$ contains all ${\mathbb P}$-null sets in ${\mathbb F}$. 
Throughout we assume that all filtrations are right continuous and complete \cite{liptser-shiryayev1977}. Define ${\mathbb F}_{T} \tri  \{ {\mathbb F}_{0,t}:   t \in [0, T]\}$. \\  

\noi In our derivations we make extensive use of the following spaces. Let $L_{{\mathbb F}_T}^2([0,T],{\mathbb R}^n) \subset   L^2( \Omega \times [0,T], d{\mathbb P}\times dt,  {\mathbb R}^n) \equiv L^2([0,T], L^2(\Omega, {\mathbb R}^n)) $ denote the space of ${\mathbb F}_{T}-$adapted random processes $\{z(t): t \in [0,T]\}$   such that
\bes
{\mathbb  E}\int_{[0,T]} |z(t)|_{{\mathbb R}^n}^2 dt < \infty,
\ees   
 which is a sub-Hilbert space of     $L^2([0,T], L^2(\Omega, {\mathbb R}^n))$.     
  Similarly, let $L_{{\mathbb F}_T}^2([0,T],  {\cal L}({\mathbb R}^m,{\mathbb R}^n)) \subset L^2([0,T] , L^2(\Omega, {\cal L}({\mathbb R}^m,{\mathbb R}^n)))$ denote the space of  ${\mathbb  F}_{T}-$adapted $n\times m$ matrix valued random processes $\{ \Sigma(t): t \in [0,T]\}$ such that
 \bes
  {\mathbb  E}\int_{[0,T]} |\Sigma(t)|_{{\cal L}({\mathbb R}^m,{\mathbb R}^n)}^2 dt  \tri  {\mathbb E} \int_{[0,T]} tr(\Sigma^*(t)\Sigma(t)) dt < \infty.
  \ees

\subsection{Distributed Stochastic System}
\label{distributed}
  
Next, we introduce the mathematical formulation of the stochastic system.    
  \noi  On the fixed probability space $\Big(\Omega,{\mathbb F},\{ {\mathbb F}_{0,t}:   t \in [0, T]\}, {\mathbb P}\Big)$ we are given a distributed stochastic dynamical decision system. It  consists of an interconnection of $N$ subsystems, and  each subsystem $i$ has,  state space ${\mathbb R}^{n_i}$, DM action space ${\mathbb A}^i \subset {\mathbb R}^{d_i}$, an exogenous  state noise space ${\mathbb W}^i \tri {\mathbb R}^{m_i}$,  an exogenous  measurement noise space ${\mathbb B}^i \tri {\mathbb R}^{k_i}$, and   initial state $x^i(0)=x_0^i$, defined by

\bi
\item[(S1)] $x^i(0)=x_0^i$:  an ${\mathbb R}^{n_i}$-valued  Random Variable;

\item[(S2)] $\{ W^i(t): t \in [0,T]\}$: an  ${\mathbb R}^{m_i}$-valued  standard Brownian motion which models the exogenous state noise, adapted to ${\mathbb F}_{T}$, independent of $x^i(0)$;

\item[(S3)] $\{ B^i(t): t \in [0,T]\}$: an  ${\mathbb R}^{k_i}$-valued  standard Brownian motion which models the exogenous measurement noise,   adapted to ${\mathbb F}_{T}$, independent of  $\{ W^i(t): t \in [0,T]\}$.

\ei
\noi The DM $\{u^i: i \in {\mathbb Z}_N\}$ take values in closed convex subsets of  metric spaces $\{({\mathbb M}^i,d): i \in {\mathbb Z}_N\}$. \\
The decentralized partial information structure available to DM $u^i$  is   generated by  noisy observation 
\bea
y^i(t) = \int_{0}^t h^i(s,x^1(s), \ldots x^N(s),y^1(s),\dots y^N(s))ds +\int_{0}^t D^{i, \frac{1}{2}}(s) dB^{i}(s), \hst  t \in [0,T], \hso \forall i \in {\mathbb Z}_N, \label{pds5}
\eea
where $x^i \in {\mathbb R}^{n_i}$ is the state of subsystem $i$ for $i =1, \ldots, N$. Notice that (\ref{pds5}) models a channel with memory and feedback. Each subsystem is described by finite dimensional coupled stochastic differential equations as follows. 
 \begin{align}
 dx^i(t) =&  f^i(t,x^i(t),u_t^i) dt  +\sigma^i(t,x^i(t),u_t^i)dW^i(t)
 + \sum_{j=1, j \neq i}^N f^{ij}(t,x^j(t),u_t^j)dt  \nonumber \\
 &+\sum_{j=1, j \neq i}^N \sigma^{ij}(t,x^j(t),u_t^j)dW^j(t) , \hst x^i(0) = x_0^i, \hso  t \in (0,T], \hso \forall i \in {\mathbb Z}_N. \label{eq1ds}
 \end{align}  
 For decentralized communication and filtering applications the right side of (\ref{eq1ds}) is independent of the DMs $u^i, i=1, \ldots, N$. \\
 Since we considered a strong  strong formulation, we define 
the filtration ${\mathbb F}_{T} \tri \{{\mathbb F}_{0,t}: t \in [0,T]\}$  as follows. Introduce the    $\sigma$-algebras  
\bes
{\mathbb F}_{0,t}^i \tri \sigma\Big\{ (x^i(0), W^i(s), B^i(s)): 0 \leq s \leq t \Big\}, \hso{\cal G}_{0,t}^{y^{i,u}} \tri \sigma\Big\{ y^i(s): 0 \leq s \leq t\Big\}, \hso t \in [0,T], \hso i=1, \ldots, N,
\ees
and the minimum $\sigma-$ algebras generated by these as follows 
\bes
 {\mathbb F}_{0,t} \tri \bigvee_{i=1}^N {\mathbb F}_{0,t}^i,  \hso  {\cal G}_{0,t}^{y^u} \tri \bigvee_{i=1}^N {\cal G}_{0,t}^{y^{i,u}}, \hso t \in [0,T].
\ees 
Next, we introduce the admissible sets of decentralized decision strategies considered in this paper. 

\bi

\item[\bf (FIS):] {\bf Feedback Information Structures.} Let  ${\cal G}_{T}^{y^{i,u}} \tri \{  {\cal G}_{0,t}^{y^{i,u}}: t \in [0, T]\} \subset \{{\mathbb F}_{0,t}: t \in [0,T]\}$ denote the information available to DM $i$, $ \forall i \in {\mathbb Z}_N$.   
     The admissible set of decentralized feedback strategies for DM $i$ is  defined by
 \begin{align}
 {\mathbb U}^{y^{i,u}}[0, T] \tri \Big\{   u^i   \in  L_{{\cal G}_T^{y^{i,u}}}^2([0,T],{\mathbb R}^{d_i})  : \:   u_t^i \in {\mathbb A}^i \subset {\mathbb R}^{d_i}, \: a.e. t \in [0,T], \: {\mathbb P}-a.s. \Big\}, \hso  \forall  i \in {\mathbb Z}_N, \label{cs1a}
     \end{align}
 where   ${\mathbb U}^{y^{i,u}}[0, T]$ is a  closed convex subset of  $L_{{\mathbb F}_T}^2([0,T],{\mathbb R}^n)$, for $i=1,2, \ldots, N$. \\
  Thus,    an $N$ tuple of DM strategies    is by definition $$(u^1,u^2, \ldots, u^N) \in {\mathbb U}^{(N),y^u}[0,T] \tri \times_{i=1}^N {\mathbb U}^{y^{i,u}}[0,T],$$ and hence it is a family of $N$ functions, say, $\Big(\mu_t^1(\cdot), \mu_t^2(\cdot), \ldots, \mu_t^N(\cdot)\Big), t \in [0,T]$,  which are nonanticipative  with respect to the information structures $\{  {\cal G}_{0,t}^{y^{i,u}}: t \in [0, T]\} , i=1,2, \ldots, N$.  
  
  The information structure of each DM is decentralized, and may be generated by local or global subsystem observables.

\item[\bf (IIS):] {\bf  Innovations Information Structures.} Let  ${\cal G}_T^{I^{i,u}} \tri  \sigma\{ I^{i}(t): 0 \leq t \leq T\}$ denote the information available to DM $i, \forall i \in {\mathbb Z}_N$, where   $\{I^i(t): t \in [0,T]\}$ is the innovations  of the process $\{y^i(t): t \in [0,T]\}$ defined by
 \begin{align}
I^{i}(t) \tri y^i(t) - \int_{0}^t {\mathbb E}\Big\{ h^{i}(s,x^1(s),\ldots, x^N(s), y^1(s), \ldots, y^N(s)) | {\cal G}_{0,s}^{y^{i,u}}  \Big\}ds, \: t \in (0,T], \: \forall i \in {\mathbb Z}_N, \label{eqsi1}
 \end{align}  
The admissible set of decentralized innovations strategies for DM $i$  is defined by 
\bea
{\mathbb U}^{I^{i,u}}[0,T] \tri \Big\{    u^i \in  {\mathbb U}^{y^{i,u}}[0, T]   : u_t^i \hso \mbox{is} \hso {\cal G}_{0,t}^{I^{i,u}}-\mbox{adapted a.e.} t \in [0,T], {\mathbb P}-a.s. \Big\}. \label{fstrategies3}
\eea
 An $N$ tuple of DM strategies    is by definition $(u^1,\ldots, u^N) \in {\mathbb U}^{(N),I^u}[0,T] \tri \times_{i=1}^N {\mathbb U}^{I^{i,u}}[0,T]$.

\ei

Define the augmented vectors by
\bes
 W \tri (W^1, \ldots, W^{N}) \in {\mathbb R}^m, \: B \tri (B^1, \ldots, B^{N}) \in {\mathbb R}^k, \: u\tri (u^1,  \ldots, u^{N}) \in {\mathbb R}^d, \: x \tri (x^1, \ldots, x^{N}) \in {\mathbb R}^n.
 \ees
  The distributed   stochastic system dynamics  are described in compact form by 
 \begin{eqnarray}
 dx(t) =  f(t,x(t),u_t) dt + \sigma(t,x(t),u_t)~dW(t), \hst x(0) = x_0, \hst  t \in (0,T], \label{eq1}
 \end{eqnarray}  
 where $f: [0,T] \times {\mathbb R}^n\times {\mathbb A}^{(N)} \longrightarrow {\mathbb R}^n$ denotes the drift and $\sigma : [0,T] \times {\mathbb R}^n\times {\mathbb  A}^{(N)} \longrightarrow {\cal L}({\mathbb R}^m, {\mathbb R}^n)$ the diffusion coefficients. \\
 The distributed observation equations  are described by the observation equation
 \bea
y(t) = \int_{0}^t h(s,x(s),y(s))ds +\int_{0}^t D^{ \frac{1}{2}}(s) dB(s), \hst  t \in [0,T], \label{eq1b}
\eea
where $h: [0,T] \times {\mathbb R}^n\times {\mathbb R}^k \longrightarrow {\mathbb R}^k$ is a function of the observations $\{y(t): 0 \leq t \leq T\}$. 

\subsection{Pay-off Functional and Team Games}
\label{tg}
   
\noi   Consider the distributed system (\ref{eq1}), (\ref{eq1b}) with decentralized partial information structures. 
  \noi Given a $u \in {\mathbb U}^{(N),y^u}[0,T]$, define the reward or performance criterion by 
 \begin{align} 
 J(u) \equiv  J(u^1,u^2,\ldots,u^N)   \tri
    {\mathbb E} \biggl\{   \int_{0}^{T}  \ell(t,x(t),u_t) dt + \varphi(x(T)\biggr\},            \label{cfd}
  \end{align}
  where $\ell: [0,T] \times {\mathbb R}^n\times {\mathbb U}^{(N)} \longrightarrow (-\infty, \infty]$ denotes the running cost function  and $\varphi : {\mathbb R}^n \longrightarrow (-\infty, \infty]$,  the terminal cost function.   
 Notice that the performance of the strategies is graded by  a single pay-off functional.

\noi The  distributed stochastic team optimization    problem with $N$ DM   is defined below.

 \begin{problem}(Team Optimality)
 \label{problemdet1}
Given the  pay-off functional (\ref{cfd}),   constraints (\ref{eq1}), (\ref{eq1b})  the  $N$ tuple of strategies   $u^o \tri (u^{1,o}, u^{2,o}, \ldots, u^{N,o}) \in {\mathbb U}^{(N),y^u}[0,T]$  is called team optimal if it satisfies  
 \bea
 J(u^{1,o}, u^{2,o}, \ldots, u^{N,o}) \leq J(u^1, u^2, \ldots, u^N),  \hst \forall  u\tri (u^1, u^2, \ldots, u^N) \in {\mathbb U}^{(N),y^u}[0,T] \label{cfd1a}
 \eea
 Any $u^o   \in {\mathbb U}^{(N),y^u}[0,T]$ satisfying (\ref{cfd1a})
is called an optimal  decision strategy (or control) and the corresponding $x^o(\cdot)\equiv x(\cdot; u^o(\cdot)), y^o(\cdot)\equiv y(\cdot; u^o(\cdot))$ (satisfying (\ref{eq1}), (\ref{eq1b})) are called an optimal state process and observation process, respectively.\\
Similarly, for $u^o \tri (u^{1,o}, u^{2,o}, \ldots, u^{N,o}) \in {\mathbb U}^{(N),I^u}[0,T]$.
  \end{problem}
 
 \noi By definition, Problem~\ref{problemdet1}  is a dynamic team problem with each DM having a different information structure (decentralized).
     An alternative approach to handle such problems with decentralized information structures is to restrict the definition of optimality to   the so-called person-by-person equilibrium.  \\
 Define 
 \bes
 \tilde{J}(v,u^{-i})\tri J(u^1,u^2,\ldots, u^{i-1},v,u^{i+1},\ldots,u^N), \hst \forall \in \in {\mathbb Z}_N.
 \ees
 
 \begin{problem}(Person-by-Person Optimality)
 \label{problemdet2}
Given the   pay-off functional (\ref{cfd}),   constraints (\ref{eq1}), (\ref{eq1b})  the  $N$ tuple of strategies   $u^o \tri (u^{1,o}, u^{2,o}, \ldots, u^{N,o}) \in {\mathbb U}^{(N),y^u}[0,T]$  is called person-by-person optimal  if it satisfies
\begin{align}
 \tilde{J}(u^{i,o}, u^{-i,o}) \leq \tilde{J}(u^{i}, u^{-i,o}), \hst \forall u^i \in {\mathbb  U}^{y^{i,u}}[0,T], \hso \forall i \in {\mathbb Z}_N. \label{cfd2}
 \end{align}
Similarly for $u^o \tri (u^{1,o}, u^{2,o}, \ldots, u^{N,o}) \in {\mathbb U}^{(N),I^u}[0,T]$.
  \end{problem}
 
 \noi The interpretaion of (\ref{cfd2}) is that the  variation of the $i$-th player is done while the rest of the players assume their optimal strategies.\\

In the next remark, the previous team games formulation is discussed in the context of distributed estimation.

\begin{remark}
\label{remf}
In distributed estimation  each subsystem is described by  unforced coupled stochastic differential equations 
 \begin{align}
 dx^i(t) =&  f^i(t,x^i(t)) dt  +\sigma^i(t,x^i(t))dW^i(t)
 + \sum_{j=1, j \neq i}^N f^{ij}(t,x^j(t))dt  \nonumber \\
 &+\sum_{j=1, j \neq i}^N \sigma^{ij}(t,x^j(t))dW^j(t) , \hst x^i(0) = x_0^i, \hso  t \in (0,T], \hso \forall i \in {\mathbb Z}_N, \label{eq1dsf}
 \end{align}  
while the observations for each subsystem are described by (\ref{pds5}). The distributed estimation objective is to determine an $N$ tuple of decision strategies   $u^o \tri (u^{1,o}, u^{2,o}, \ldots, u^{N,o}) \in {\mathbb U}^{(N),y^u}[0,T]$ which is team optimal or person-by-person optimal (according to Problems~\ref{problemdet1}, \ref{problemdet2}) subject to constraints (\ref{eq1dsf}), (\ref{pds5}). This distributed estimation problem formulated via team theory,  is a generalization of the static team theory discussed in \cite{radner1962,ho1980,krainak-speyer-marcus1982a,krainak-speyer-marcus1982b}. However, we point out that for  distributed filtering there is no reason to consider innovations information structures.  
\end{remark}

   \section{Strong Solutions and Variational Equation}
  \label{sexistence}
In this section we introduce assumptions which will allow us to show existence of strong ${\mathbb F}_{T}-$adapted continuous solutions to (\ref{eq1}) and (\ref{eq1b}). We also introduce the variational equation which is utilized to derive the stochastic minimum principle using the methodology in \cite{ahmed-charalambous2012a,charalambous-ahmedFIS_Parti2012}. 
 
\noi Let $B_{{\mathbb F}_T}^{\infty}([0,T], L^2(\Omega,{\mathbb R}^n))$ denote the space of ${\mathbb F}_{T}$-adapted ${\mathbb R}^n$ valued second order random processes endowed with the norm topology  $\parallel  \cdot \parallel$ defined by  
\bes
 \parallel x\parallel^2  \tri \sup_{t \in [0,T]}  {\mathbb E}|x(t)|_{{\mathbb R}^n}^2.
 \ees
The existence of strong solution is based on the following assumptions.

\begin{assumptions}(Main assumptions)
\label{A1-A4}
The coefficients of the state and observation equations (\ref{eq1}), (\ref{eq1b}) 
are     Borel measurable  maps: 
\bes
 f: [0,T] \times {\mathbb R}^n \times {\mathbb A}^{(N)} \longrightarrow {\mathbb R}^n , \hst  \sigma: [0,T] \times {\mathbb R}^n \times {\mathbb A}^{(N)} \longrightarrow {\cal L}({\mathbb R}^m, {\mathbb R}^n), 
\ees
\bes 
 h^i : [0,T] \times {\mathbb R}^n \times {\mathbb R}^k \longrightarrow {\mathbb R}^{k_i}, \hso \forall i \in {\mathbb Z}_N.
 \ees
These satisfy the following basic conditions. 

\noi There exists a $K  >0$ such that

\begin{description}
\item[(A1)] $|f(t,x,u)-f(t,z,u)|_{{\mathbb R}^n} \leq K |x-z|_{{\mathbb R}^n}$ uniformly in $u \in {\mathbb A}^{(N)}$;

\item[(A2)] $|f(t,x,u)-f(t,x,v)|_{{\mathbb R}^n} \leq K |u-v|_{{\mathbb R}^d}$ uniformly in $x \in {\mathbb R}^{n}$;

\item[(A3)] $|f(t,x,u)|_{{\mathbb R}^n} \leq K (1 + |x|_{{\mathbb R}^n}+|u|_{{\mathbb R}^d} )$;

\item[(A4)] $|\sigma(t,x,u)-\sigma(t,z,u)|_{{\cal L}({\mathbb R}^m, {\mathbb R}^n)} \leq K |x-z|_{{\mathbb R}^n}$ uniformly in  $u \in {\mathbb A}^{(N)}$;

\item[(A5)] $|\sigma(t,x,u)-\sigma(t,x,v)|_{{\cal L}({\mathbb R}^m, {\mathbb R}^n)} \leq K |u-v|_{{\mathbb R}^d}$ uniformly in  $x \in {\mathbb R}^{n}$;

\item[(A6)] $|\sigma(t,x,u)|_{{\cal L}({\mathbb R}^m, {\mathbb R}^n)} \leq K (1+ |x|_{{\mathbb R}^n}+|u|_{{\mathbb R}^d})$;

\item[(A7)] $|h(t,x,y)|_{{\mathbb R}^k} \leq K(1 + |x|_{{\mathbb R}^n}  +  |y|_{{\mathbb R}^k}   )$;

\item[(A8)]  $|h(t,x,y)-h(t,z,\tilde{y})|_{{\mathbb R}^k} \leq K \Big( |x-z|_{{\mathbb R}^{n}} + |y-\tilde{y}|_{{\mathbb R}^{k}} \Big)$;

\item[(A9)] For any $a, b \in C([0,T],{\mathbb R}^{k_i})$ the nonanticipative mapping $u^i \equiv \mu^i : [0,T] \times C([0,T],{\mathbb R}^{k_i}) \longrightarrow {\mathbb A}^i$ satisfies the Lipschitz condition
\bes
|\mu^i(t,a) -\mu^i(t,b)|_{{\mathbb R}^{d_i}} \leq K |a-b|_{C([0,T], {\mathbb R}^{k_i})}, \hst i=1, \ldots, N.
\ees 

\item[A10)] $D^{i,\frac{1}{2}}(t) >0, \forall t \in [0,T]$ and $D^{i,\frac{1}{2}}(\cdot)$ is uniformly bounded, $\forall i \in {\mathbb Z}_N$.

\end{description}
\end{assumptions}

\noi The following lemma proves  the  existence of solutions and their continuous dependence on the decision variables.

\begin{lemma}
\label{lemma3.1}
 Suppose Assumptions~\ref{A1-A4} hold. Then for any ${\mathbb  F}_{0,0}$-measurable initial state $x_0$ having finite second moment, and any $u \in {\mathbb U}^{(N),y^u}[0,T]$,  the following hold. 

\begin{description} 
 
 \item[(1)]  System (\ref{eq1}), (\ref{eq1b}) has a unique solution   $(x,y) \in B_{{\mathbb F}_T}^{\infty}([0,T],L^2(\Omega,{\mathbb R}^{n+k}))$  having a continuous modification, that is, $(x,y) \in C([0,T],{\mathbb R}^{n+k})$,  ${\mathbb P}-$a.s, $\forall i \in {\mathbb Z}_N$.
 
\item[(2)]  The solution of  system  (\ref{eq1}), (\ref{eq1b}) is continuously dependent on the control, in the sense that, as $u^{i, \alpha}  \longrightarrow u^{i,o}$  in ${\mathbb U}^{y^{i,u}}[0,T]$, $\forall i \in {\mathbb Z}_N$,  $(x^\alpha, y^\alpha) \longrightarrow (x^o,y^o) $ in $B_{{\mathbb F}_T}^{\infty}([0,T],L^2(\Omega,{\mathbb R}^{n+k}))$, $\forall i \in {\mathbb Z}_N$.
\end{description}
Similarly for $u \in {\mathbb U}^{(N),I^u}[0,T]$.

\end{lemma} 

\begin{proof} {\bf (1)} 
Consider the augmented system $X\tri (x,y)$ and the associated stochastic differential equation of $X$. The proof for the first part of the lemma is classical and hence omitted.

\noi {\bf (2)}  Next, we consider the second part asserting the  continuity of $u$ to solution map  $u\longrightarrow (x,y).$  Let $\{\{u^{i,\alpha}: i=1,2,\ldots, N\},u^o\}$ be any pair of  DM strategies from ${\mathbb U}^{(N),y^u}[0,T] \times {\mathbb U}^{(N),y^u}[0,T]$ and  $\{x^\alpha,y^\alpha, x^o,y^o\}$ denote the corresponding pair of solutions of the system (\ref{eq1}), (\ref{eq1b}).  Let $u^{i,\alpha} \longrightarrow u^{i,o}, i=1,2,\ldots,N$.  We must show that  $(x^\alpha, y^\alpha) \longrightarrow (x^o,y^o)$ in  $B_{{\mathbb F}_T}^{\infty}([0,T], L^2(\Omega,R^{n+k})).$    By the definition of solution to (\ref{eq1}), it can be verified that 
\begin{align} 
x^\alpha(t)-x^o(t) =& \int_0^t \Big\{f(s,x^\alpha(s),u^\alpha_s)-f(s,x^o(s),u_s^\alpha)\Big\} ds  \nonumber \\
&+ \int_0^t \Big\{\sigma(s,x^\alpha(s),u^\alpha_s)-\sigma(s,x^o(s),u_s^\alpha)\Big\} dW(s) + e_{1}^\alpha(t) + e_{2}^\alpha(t),  \hst t \in [0,T], 
 \label{eq3}
 \end{align} 
 where
\begin{align}
e_{1}^\alpha(t) =&  \int_0^t \Big\{ f(s,x^o(s),u^\alpha_s)-f(s,x^o(s),u_s^o)] ds  \label{eqe1}  \\ 
  e_{2}^\alpha(t) =& \int_0^t \Big\{ \sigma(s,x^o(s),u^\alpha_s)-\sigma(s,x^o(s),u_s^o)\Big\} dW(s). \label{eqe2}
  \end{align}
   Using  the standard martingale inequality into (\ref{eq3}),  it follows from it  and {\bf (A1), (A4)} that there exist constants $C_1,C_2> 0$ such that
\begin{eqnarray}
 {\mathbb E} |x^\alpha(t)-x^o(t)|_{{\mathbb R}^n}^2 \leq C_1 \int_0^t K^2  {\mathbb E}|x^\alpha(s)-x^o(s)|_{{\mathbb R}^n}^2 + C_2 \bigl( {\mathbb E}|e_1^\alpha(t)|_{{\mathbb R}^n}^2 + {\mathbb E}|e_2^\alpha(t)|_{{\mathbb R}^n}^2\bigr). \label{eq4} \end{eqnarray} 
 Clearly, by the Cauchy-Schwartz inequality, and martingale inequality, it follows from   {\bf (A2), (A5)} that  
\begin{align}
{\mathbb E}|e_1^\alpha(t)|_{{\mathbb R}^n}^2 \leq & T \: {\mathbb E}\int_0^t  |f(s,x^o(s),u_s^\alpha)-f(s,x^o(s),u_s^o)|_{{\mathbb R}^n}^2 ds \leq T  {\mathbb E} \int_0^t K^2  |u_s^\alpha-u_s^o)|_{{\mathbb R}^d}^2 ds, \label{ineq1}\\
  {\mathbb E}|e_2^\alpha(t)|_{{\mathbb R}^n}^2 \leq & 4 \: {\mathbb  E}\int_0^t |\sigma(s,x^o(s),u_s^\alpha)-\sigma(s,x^o(s),u_s^o)|_{{\cal L}({\mathbb R}^m,{\mathbb R}^n)}^2 ds \leq 4 \: {\mathbb  E}\int_0^t K^2 |u_s^\alpha-u_s^o|_{{\mathbb R}^d}^2 ds  . \label{ineq2}
\end{align}
Similarly, by {\bf (A8)} 
\begin{eqnarray}
 {\mathbb E} |y^\alpha(t)-y^o(t)|_{{\mathbb R}^k}^2 \leq T \int_0^t K^2  {\mathbb E} \Big(  | x^\alpha(s)-x^o(s)|_{{\mathbb R}^n}^2 + | y^\alpha(s)-y^o(s)|_{{\mathbb R}^k}^2  \Big)ds. \label{eq4a} 
 \end{eqnarray} 
The  integrands  in the right side of  inequalities (\ref{ineq1}), (\ref{ineq2}) converge to zero  for almost all $s \in [0,T], {\mathbb  P}-$a.s. Moreover, these integrands  are dominated by integrable functions.   Hence, by Lebesgue dominated convergence theorem  the terms  $\{e_1^\alpha, e_2^\alpha\}$ converge to zero uniformly on $[0,T]$. Define $\rho^\alpha(t)\tri {\mathbb E}  \Big(| x^\alpha(s)-x^o(s)|_{{\mathbb R}^n}^2 + | y^\alpha(s)-y^o(s)|_{{\mathbb R}^k}^2\Big)$. Then by Gronwall inequality applied to  $\rho^\alpha$, it can be shown that $\rho^\alpha \longrightarrow 0$ as   $u^{i, \alpha}  \longrightarrow u^{i,o}$  in ${\mathbb U}^{y^{i,u}}[0,T]$, $\forall i \in {\mathbb Z}_N$. The above analysis holds for innovations information structures. This completes the derivation.
\end{proof}

\ \ \ \
\noi Throughout the paper we  assume existence of  a minimizer $u^o \in {\mathbb U}^{(N),y^u}[0,T]$ for Problem~\ref{problemdet1}. For randomize (relaxed) strategies existence can be shown as in \cite{ahmed-charalambous2012a}.

\noi Next, we prepare to  introduce the variational equation of the augmented system $(x,y)$.

\noi Define the augmented  vectors
\bes
X \tri (x,y) \in {\mathbb R}^n \times {\mathbb R}^k, \hst B \tri (B^1,B^2, \ldots, B^N) \in {\mathbb R}^k, \hst y \tri (y^1,y^2, \ldots, y^N) \in {\mathbb R}^k,
\ees
and  the augmented drift and diffusion coefficients, and terms in the pay-off associated with them by
\begin{align}
&F(t,X, u) \tri \left[ \begin{array}{c} f(t,x,u) \\ h(t,x,y)   \end{array} \right], \hso G(t,X, u) \tri \left[ \begin{array}{cc} \sigma(t,x,u) & 0 \\ 0 & D^{\frac{1}{2}}(t) \end{array} \right] , \nonumber  \\
 &h(t,x,y)  \tri\left[ \begin{array}{c} h^{1}(t,x,y) \\  \ldots \\ h^{N}(t,x,y) \end{array} \right],  \hso D^{\frac{1}{2}}(t) \tri diag \{D^{1,\frac{1}{2}}(t), \ldots, D^{N,\frac{1}{2}}(t) \}   \nonumber \\
 &L(t,X,u) \tri \ell(t,x,u), \hst \Phi(X) \tri \varphi(x). \nonumber 
\end{align}
Then  the augmented system  is expressed in compact form by
\bea
dX(t) = F(t,X(t),u_t)dt + G(t,X(t),u_t) \left[ \begin{array}{c} dW(t) \\ dB(t) \end{array} \right], \hst X(0)=X_0, \hst t \in (0,T]. \label{pds20}  
 \eea
For strategies ${\mathbb U}^{(N), y^{u}}[0, T]$,  since the state of the augmented system is $X=(x,y)$,   when considering variations of the state trajectory $X$, due to  variation of $u$,  there will be derivatives of $u(\cdot)$ with respect to $y$. To avoid this technicality we introduce the following assumptions.

\begin{assumptions}
\label{a-nf}
The diffusion coefficients $\sigma$  is restricted to   the    Borel measurable  map $\sigma: [0,T] \times {\mathbb R}^n \times {\mathbb A}^{(N)} \longrightarrow {\cal L}({\mathbb R}^n, {\mathbb R}^n)$ 
(e.g., it is independent of $u$) and 
\begin{description}
\item[(A11)] $\sigma(\cdot,\cdot)$ and $\sigma^{-1}(\cdot, \cdot)$  are bounded.
\end{description}
\end{assumptions} 
\ \

\noi Under the additional Assumptions~\ref{a-nf} we can show the following Lemma.

\begin{lemma}
\label{lemma-nf}
Consider Problem~\ref{problemdet1} under Assumptions~\ref{A1-A4}, \ref{a-nf} hold.  Define the $\sigma-$algebras
\bes 
 {\cal F}_{0,t}^{x(0),W,B} \tri \sigma\Big\{(x(0), W(s),B(s)): 0 \leq s \leq t \Big\}, \hso {\cal F}_{0,t}^{x^u,y^u} \tri \sigma\Big\{(x(s), y(s)): 0 \leq s \leq t \Big\}, \: \forall t \in [0,T].
\ees
 If  $u \in {\mathbb U}_{reg}^{(N),y^{u}}[0,T]$  then $ {\cal F}_{0,t}^{x(0),W,B}  = {\cal F}_{0,t}^{x^u,y^u}, \forall t \in [0,T]$.
\end{lemma}
\begin{proof} This follows directly from Assumptions~\ref{a-nf} and the invertibility of $D^i(\cdot), \forall i \in {\mathbb Z}_N$.
\end{proof}

\noi   Recall that $\{x(t), y(t): t \in [0,T]\}$ are the strong ${\mathbb F}_T-$adapted solutions of the state and observation equations.  Under the conditions of Lemma~\ref{lemma-nf}  for any $u^i \in {\mathbb U}^{y^{i,u}}[0,T]$ which is ${\cal G}_T^{y^{i,u}}-$adapted there exists a function $\phi^i(\cdot)$ measurable with respect  to a sub-$\sigma-$algebra of ${\cal F}_{0,t} \subset {\cal F}_{0,t}^{x(0),W,B}$ such that $u_t^i( \omega)= \phi^i(t, x(0), B(\cdot \bigwedge t,\omega), W(\cdot \bigwedge t,\omega)), {\mathbb P}-a.s.\: \omega \in {\Omega}, \forall t \in [0,T], i=1,\ldots N$.  \\
Define all such adapted nonanticipative functions by 
\bea
\overline{\mathbb U}^{i}_{na}[0, T] \tri \Big\{   u^i   \in  L_{{\cal F}_{T}^{ }}^2([0,T],{\mathbb R}^{d_i})  : \:   u_t^i \in {\mathbb A}^i \subset {\mathbb R}^{d_i}, \: a.e. t \in [0,T], \: {\mathbb P}-a.s. \Big\}, \hso \forall  i \in {\mathbb Z}_N. \label{na1}
     \eea

\noi Next, we introduce the following additional assumptions.

\begin{assumptions}
\label{a-nf1r}
 ${\mathbb U}^{y^{i,u}}[0, T]$ is dense in  $\overline{\mathbb U}_{na}^{i}[0, T], \forall i \in {\mathbb Z}_N$.
\end{assumptions}

\noi Under Assumptions~\ref{a-nf1r} we can show the following theorem.

\begin{theorem}
\label{thm-nf}
Consider Problem~\ref{problemdet1} under Assumptions~\ref{A1-A4}, \ref{a-nf1r}. Further, assume $\ell$ is Borel measurable, continuously differentiable with respect to  $(x,u)$, and $\varphi$ is continously differentiable with respect to $x$, and there exist $K_1, K_2 >0$ such that
\bes
|\ell_x(t,x,u)|_{{\mathbb R}^n}+|\ell_u(t,x,u) |_{{\mathbb R}^d}     |  \leq K_1 \big(1+|x|_{{\mathbb R}^n} + |u|_{{\mathbb R}^d} \big),  \hso |\varphi_x(T,x)|_{\mathbb R} \leq K_2 \big(1+ |x|_{{\mathbb R}^n}\big).
\ees
Then
\bea 
 \inf_{ u \in  \times_{i=1}^N \overline{\mathbb U}_{na}^{i}[0,T]} J(u)=  \inf_{ u  \in \times_{i=1}^N  {\mathbb U}^{y^{i,u}}[0,T]} J(u). \label{eqcost}
 \eea
\end{theorem}

\begin{proof} Since Assumptions~\ref{a-nf1r} holds, it is sufficient to show that  as $u^{i, \alpha} \longrightarrow u^{i}$  in $\overline{\mathbb U}_{na}^i[0,T]$, $\forall i \in {\mathbb Z}_N$,  then $J(u^\alpha) \longrightarrow J(u)$. From the derivation of Lemma~\ref{lemma3.1}, we can  show that ${\mathbb E} \sup_{s \in [0,t] } |x^\alpha(s)-x(s)|_{{\mathbb R}^n}$ converges to zero as $\alpha \longrightarrow \infty$, hence it is sufficient to show that $|J(u^\alpha)-J(u)|$ also converges to zero, as $\alpha \longrightarrow \infty$. By the assumptions on $\{\ell, \varphi\}$, and by the mean value theorem we have the following inequality.
\begin{align}
| J(u^\alpha)-J(u)| \leq & K_1 \: {\mathbb E} \Big\{ \int_{[0,T]}^{} \Big( |x^\alpha(t)|_{{\mathbb R}^n} + |u_t^\alpha|_{{\mathbb R}^d} +  |x(t)|_{{\mathbb R}^n} + |u_t|_{{\mathbb R}^d} +1 \Big) \nonumber \\
&.\Big(  |x^\alpha(t)-x(t)|_{{\mathbb R}^n} + |u_t^\alpha-u_t|_{{\mathbb R}^d} \Big)dt\Big\} \nonumber \\
&+ K_2 {\mathbb E} \Big\{ \Big( |x^\alpha(T)|_{{\mathbb R}^n} + |x(T)|_{{\mathbb R}^n} +1 \Big) |x^\alpha(T)-x(t)|_{{\mathbb R}^n} \Big\}. \label{nac}
\end{align}
Since ${\mathbb E} \sup_{s \in [0,t] } |x^\alpha(s)-x(s)|_{{\mathbb R}^n} \longrightarrow 0$ as $\alpha \longrightarrow \infty$, then $|J(u^\alpha)-J(u)|$ also converges to zero, as $\alpha \longrightarrow \infty$.

\end{proof}

\noi The point to be made regarding Theorem~\ref{thm-nf} is that if $u \in {\mathbb U}^{(N),y^{u}}[0,T]$ achieves the infimum of $J(u)$ then it is also optimal with respect to some measurable functionals of subsets  $\{(x(0), (W(s), B(s)): 0 \leq s \leq T\}$. Consequently, the  necessary conditions for $u \in {\mathbb U}^{(N),y^{u}}[0,T]$ to be optimal are those for which $u \in \times_{i=1}^N \overline{{\mathbb U}}_{na}^{i}[0,T]$ is optimal. 

\begin{remark}
\label{innovations}
Strategies adapted to the innovations process $u^i \in {\cal G}_T^{I^{i,u}}$, $\forall i \in {\mathbb Z}_N$  are often utilized  to derive the separation theorem of partially observed stochastic control problems. It is well known that $\{I^{i,u}(t): t \in [0,T]\}$ is an ${\cal G}_T^{y^{i,u}}-$adapted Wiener process. 
    Moreover, if the innovations process and observation process generate the same $\sigma-$algebra,  ${\cal G}_{0,s}^{y^{i,u}}= {\cal G}_{0,s}^{I^{i,u}}$, and the innovations process  is independent of $u$, then $I^{i,u}(t)=I^{i,0}(t),  \forall t \in [0,T]$. Define the $\sigma-$algebra ${\cal G}_{0,t}^{I^{i,0}} \tri \sigma\{ I^{i,0}(s): 0 \leq s \leq t \}, t \in [0,T]$. Under these conditions the  necessary conditions for $u \in {\mathbb U}^{(N),y^{u}}[0,T]$ to be optimal are those for which  $u \in {\mathbb U}^{(N),I^{0}}[0,T]$, defined by 
\bea
{\mathbb U}^{(N),I^{0}}[0,T] \tri \Big\{ {\mathbb U}^{(N),y^u}: u_t^i\: \: \mbox{is} \:  \: {\cal G}_{0,t}^{I^{i,0}}-\mbox{ adapted},  \: \forall i \in {\mathbb Z}_N\Big\}. \label{sin}
\eea
Note that Assumptions~\ref{a-nf1r} are not required  for distributed filtering applications because   the decentralized information structures $ {\cal G}_T^{y^{i}}$ are independent of $u$, and hence it is not very difficult to show $ {\cal G}_T^{y^{i}}={\cal G}_{T}^{I^{i,0}}, \forall i \in {\mathbb Z}_N$.   \\
After deriving the necessary conditions we also show that under certain convexity conditions that these are also sufficient. Consequently, for the sufficient part we do not require Assumptions~\ref{a-nf1r}.
\end{remark}

\noi  For the  derivation of stochastic minimum principle of  optimality   we shall require stronger regularity conditions on the maps  $\{ f,\sigma, h\}$, as well as, for the running and terminal pay-offs functions   $\{\ell,\varphi\}.$  These  are given below.

\begin{assumptions}
\label{NCD1}  
${\mathbb E}|x(0)|_{{\mathbb R}^n} < \infty$ and the maps of $\{f,\sigma,\ell, \varphi\} $ satisfy the following conditions.

\begin{description}

\item[(B1)] The map $f: [0,T] \times {\mathbb R}^n \times {\mathbb A}^{(N)} \longrightarrow {\mathbb R}^n$ is continuous in $(t,x,u)$ and continously differentiable with respect to $(x,u)$;

\item[(B2)] The map $\sigma: [0,T] \times {\mathbb R}^n \times {\mathbb A}^{(N)} \longrightarrow {\cal L}({\mathbb R}^m; {\mathbb R}^n)$ is continuous in $(t,x,u)$ and continously differentiable with respect to $(x,u)$;

\item[(B3)] The map $h: [0,T] \times {\mathbb R}^n \times {\mathbb R}^k \longrightarrow {\mathbb R}^k$ is continuous in $(t,x,y)$ and continously differentiable with respect to $(x,y)$;

\item[(B4)] The first derivatives $\{f_x, f_u,\sigma_x, \sigma_u \}$ are bounded  uniformly on $[0,T] \times {\mathbb R}^n \times {\mathbb A}^{(N)}$;

\item[(B5)] The first derivative $\{h_x,h_y\}$ are bounded  uniformly on $[0,T] \times {\mathbb R}^n \times  {\mathbb R}^k$;

\item[(B6)] The maps $\ell: [0,T] \times {\mathbb R}^n \times {\mathbb A}^{(N)} \longrightarrow (-\infty, \infty]$ is Borel measurable, continuously differentiable with respect to  $(x,u)$, the map $\varphi: [0,T] \times {\mathbb R}^n \longrightarrow (-\infty, \infty]$ is continously differentiable with respect to $x$, and there exist $K_1, K_2 >0$ such that
\bes
|\ell_x(t,x,u)|_{{\mathbb R}^n}+|\ell_u(t,x,u) |_{{\mathbb R}^d}     |  \leq K_1 \big(1+|x|_{{\mathbb R}^n} + |u|_{{\mathbb R}^d} \big),  \hso |\varphi_x(T,x)|_{{\mathbb R}^n}   \leq K_2 \big(1+ |x|_{{\mathbb R}^n}\big)
\ees

\item[(B7)] Conditions {\bf   (A9), (A10)} of Assumptions~\ref{A1-A4} hold.

\end{description}

\end{assumptions}

\noi Consider the Gateaux derivative of $G$ with respect to the  variable  at the point $(t,z,v) \in [0,T] \times {\mathbb R}^{n+k}\times {\mathbb A}
^{(N)}$   in the direction $\eta \in {\mathbb R}^{n+k}$ defined by
\bes
   G_X(t,z,v; \eta) \tri   \lim_{\varepsilon \rightarrow 0}\frac{1}{\varepsilon} \Big\{ G(t,z + \varepsilon \eta, \nu)- G(t,z,v)\Big\}, \hst t \in [0,T].
  \ees 
    Note that the map $\eta \longrightarrow G_X(t,z,\nu; \eta)$ is linear, and  it follows from Assumptions~ \ref{NCD1}, {\bf (B3), (B5)} that  there exists a finite positive number $\beta>0$ such that
    \bes
     |G_X(t,z,\nu; \eta)|_{{\cal L}({\mathbb R}^{m+k},{\mathbb R}^{n+k})} \leq \beta |\eta|_{{\mathbb R}^{n+k}},   \hst t \in  [0,T].
     \ees
\noi   In order to present the necessary conditions of optimality we need the so called variational equation.
Suppose $u^o \tri (u^{1,o}, u^{2,o}, \ldots, u^{N,o}) \in {\mathbb U}^{(N),I^u}[0,T]$ denotes the optimal decision and $u \tri (u^1, u^2, \ldots, u^n) \in {\mathbb U}^{(N),I^u}[0,T]$ any other decision.  Since ${\mathbb U}^{I^{i,u}}[0,T]$ is convex $\forall i \in {\mathbb Z}_N$, it is clear that  for any $\varepsilon \in [0,1]$, 
\bes
 u_t^{i,\varepsilon} \tri u_t^{i,o} + \varepsilon (u_t^i-u_t^{i,o}) \in {\mathbb U}^{I^{i,u}}[0,T], \hst \forall i \in {\mathbb Z}_N.
 \ees
 Let $X^{\varepsilon}(\cdot)\equiv X^\veps(\cdot; u^\veps(\cdot))$ and  $X^{o}(\cdot) \equiv X^o(\cdot;u^o(\cdot))  \in B_{{\mathbb F}_T}^{\infty}([0,T],L^2(\Omega,{\mathbb R}^{n+k}))$ denote the solutions  of the system equation (\ref{pds20})  corresponding to  $u^{\varepsilon}(\cdot)$ and $u^o(\cdot)$, respectively.  Consider the limit 
 \bes
  Z(t) \tri \lim_{\varepsilon\downarrow 0}  \frac{1}{\veps} \Big\{X^{\varepsilon}(t)-X^o(t)\Big\} , \hst t \in [0,T]. 
  \ees
We have the following result characterizing the the variational process $\{Z(t): t \in [0,T]\}$.

\begin{lemma}
\label {lemma4.1}
Suppose  Assumptions~\ref{NCD1} hold. For strategies ${\mathbb U}^{(N),I^u}[0,T]$ the process $\{Z(t): t \in [0,T]\}$ is an element of the Banach space    $B_{{\mathbb F}_T}^{\infty}([0,T],L^2(\Omega,{\mathbb R}^{n+k}))$ and it  is the unique solution of the variational stochastic differential equation
 \begin{align} 
 &dZ(t) = F_X(t,X^o(t),u_t^o)Z(t)dt + G_X(t,X^o(t),u_t^o; Z(t)) \left[ \begin{array}{c} dW(t)\\dB(t) \end{array} \right] \nonumber \\ 
 &+ \sum_{i=1}^N F_{u^i}(t,X^o(t),u_t^{,o})(u_t^i-u_t^{i,o})dt + \sum_{i=1}^N G_{u^i}(t,X^o(t),u_t^{o};u_t^i-u_t^{i,o}) \left[\begin{array}{c}   dW(t) \\ dB(t) \end{array} \right], \hst Z(0)=0. \label{eq9}   
 \end{align} 
 having a continuous modification.\\  
Under the addition Assumptions~\ref{a-nf1r}  the above statements hold for strategies ${ \mathbb U}^{(N),y^{u}}[0, T]$. \\
Moreover, (\ref{eq9}) is the variational equation  for distributed filtering applications (without imposing Assumptions~\ref{a-nf1r}).
  \end{lemma}

\begin{proof}  This follows from  \cite{charalambous-ahmedFIS_Parti2012} by considering the augmented system.

  \end{proof}

\noi Using the variation equation of Lemma~\ref{lemma4.1}, we note that the results given in \cite{charalambous-ahmedFIS_Parti2012} for nonrandomized strategies are directly applicable to the augmented system (\ref{pds20}). In fact one can also consider  randomized strategies as in \cite{charalambous-ahmedFIS_Parti2012}.

\section{Optimality Conditions for Noisy Information Structures }
\label{smp}
In this section we derive necessary and sufficient optimality conditions for the team game of Problem~\ref{problemdet1}. In view of the results obtained in the previous section, specifically, Lemma~\ref{lemma4.1}, the  stochastic  minimum principle of  optimality  for Problem~\ref{problemdet1}, described in terms of the augmented system (\ref{pds20}),  follows directly from  the results in \cite{charalambous-ahmedFIS_Parti2012}.\\
 
  Before we introduce the optimality conditions we define the Hamiltonian system of equations. To this end, define the Hamiltonian 
\bes
 {\mathbb  H}: [0, T] \times {\mathbb R}^n\times {\mathbb R}^n\times {\cal L}({\mathbb R}^m,{\mathbb R}^n)\times  {\mathbb A}^{(N)} \longrightarrow {\mathbb R}
\ees  
   by  
   \begin{align}
    {\mathbb H} (t,x,\psi,q_{11},u) \tri    \langle f(t,x,u),\psi \rangle + tr (q_{11}^*\sigma(t,x))
     + \ell(t,x,u),  \hst  t \in  [0, T]. \label{h1}
    \end{align}
    \noi For any $u \in {\mathbb U}^{(N),y^u}[0,T],   {\mathbb U}^{(N),I^u}[0,T],$ the adjoint process is $$(\psi,q_{11},q_{12}) \in  L_{{\mathbb F}_T}^2([0,T], {\mathbb R}^n)\times L_{{\mathbb F}_T}^2([0,T] ,{\cal L}({\mathbb R}^m,{\mathbb R}^n)) \times L_{{\mathbb F}_T}^2([0,T] ,{\cal L}({\mathbb R}^k,{\mathbb R}^n)) $$ and satisfies the following backward stochastic differential equation
\begin{align} 
d\psi (t)  =& -f_x^{*}(t,x(t),u_t)\psi (t)  dt - V_{q_{11}}(t) dt -\ell_x(t,x(t),u_t) dt \nonumber  \\
 &+ q_{11}(t) dW(t) + q_{12}(t)dB(t), \hst t \in [0,T),    \nonumber    \\ 
=&- {\mathbb H}_x (t,x(t),\psi(t),q_{11}(t),u_t) dt + q_{11}(t)dW(t) + q_{12}(t)dB(t),    \hso  t \in [0,T), \label{adj1a} \\ 
 \psi(T) =&   \varphi_x(x(T)),  \label{eq18}  
 \end{align}
  where  $V_{q_{11}} \in L_{{\mathbb F}_T}^2([0,T],{\mathbb R}^n)$ is   given by  $\langle V_{q_{11}}(t),\zeta\rangle = tr (q_{11}^*(t)\sigma_x(t,x(t); \zeta)), t \in [0,T]$ (e.g., $V_{q_{11}}(t) = \sum_{k=1}^n \Big(\sigma_x^{(k)}(t,x(t))\Big)^*q_{11}^{(k)}(t), \hst t \in [0,T],$ $\sigma^{(k)}$ is the $kth$ column of $\sigma$, $\sigma_x^{(k)}$ is the  derivative of $\sigma^{(k)}$ with respect to the state, $q_{11}^{(k)}$ is the $kth$ column of $q_{11}$, for $k=1, 2, \ldots, m$).\\ 
The state process satisfies the stochastic differential equation
\begin{align}
dx(t) &=f(t,x(t),u_t)dt + \sigma(t,x(t))dW(t), \hst t \in (0, T],  \nonumber  \\
        & = {\mathbb H}_\psi (t,x(t),\psi(t),q_{11}(t),u_t)     dt + \sigma(t,x(t)) dW(t), \hst t \in (0, T],  \label{st1a} \\
        x(0) &=  x_0 \label{st1i}  
 \end{align}
The above Hamiltonian system of equations is expressed in terms of the original distributed system of equations (\ref{eq1}), (\ref{eq1b}), and it is obtained by first deriving the Hamiltonian system of equations for the augmented system (\ref{pds20}) (we shall clarify this step in the next section).

\subsection{Necessary Conditions of Optimality}
\label{necessary}
We now prepare to derive  the necessary conditions for team optimality. Specifically, given that $u^o \in {\mathbb U}^{(N),y^u}[0,T]$ or $u^o \in {\mathbb U}^{(N),I^u}[0,T]$ is team optimal the question we address is whether it satisfies the Hamiltonian system of equations (\ref{h1})-(\ref{st1i}).

\noi By utilizing \cite{charalambous-ahmedFIS_Parti2012} we have  following  necessary conditions.

 \begin{theorem} (Necessary conditions for team optimality)
 \label{theorem5.1o}
Consider Problem~\ref{problemdet1} under Assumptions~\ref{NCD1}, and ${\mathbb A}^i$ a closed, bounded and convex subset of ${\mathbb R}^{k_i}, i=1, \ldots N$.  \\

\noi  For  an element $ u^o \in {\mathbb U}^{(N),I^u}[0,T]$ with the corresponding solution $x^o \in B_{{\mathbb F}_T}^{\infty}([0,T], L^2(\Omega,{\mathbb R}^{n}))$ to be team optimal, it is necessary  that 
the following hold.

   \begin{description}

\item[(1)]  There exists a square integrable semi martingale  $m^o$ with the intensity process $({\psi}^o,q_{11}^o,q_{12}) \in  L_{{\mathbb F}_T}^2([0,T],{\mathbb R}^{n})\times L_{{\mathbb F}_T}^2([0,T],{\cal L}({\mathbb R}^{m},{\mathbb R}^{n}))\times L_{{\mathbb F}_T}^2([0,T],{\cal L}({\mathbb R}^{k},{\mathbb R}^{n}))$.
 
 \item[(2) ]  The variational inequality is satisfied:

\begin{align}    
  \sum_{i=1}^N {\mathbb E} \int_{0}^{T}  \la {\mathbb H}_{u^i}(t,x^o(t),\psi^o(t),q_{11}^o(t), u_t^{o}),u_t^i-u_t^{i,o} \ra dt   
  \geq 0, \hst \forall u \in {\mathbb U}^{(N)^{I^u}}[0,T]. \label{eq16o}
\end{align}

\item[(3)]  The process $({\psi}^o,q_{11}^o,q_{12}^o) \in  L_{{\mathbb F}_T}^2([0,T],{\mathbb R}^{n})\times L_{{\mathbb F}_T}^2([0,T],{\cal L}({\mathbb R}^{m},{\mathbb R}^{n})) \times L_{{\mathbb F}_T}^2([0,T],{\cal L}({\mathbb R}^{k},{\mathbb R}^{n}))$ is a unique solution of the backward stochastic differential equation (\ref{adj1a}), (\ref{eq18}) such that $u^o \in {\mathbb U}^{(N),I^u}[0,T]$ satisfies  the  point wise almost sure inequalities with respect to the $\sigma$-algebras ${\cal G}_{0,t}^{I^{i,u}}   \subset {\mathbb F}_{0,t}$, $ t\in [0,T], i=1, 2, \ldots, N:$ 

\begin{align} 
&  \la {\mathbb E} \Big\{   {\mathbb H}_{u^i}(t,x^o(t),\psi^o(t),q_{11}^o(t),u_t^{o}) |{\cal G}_{0, t}^{I^{i,u^o}} \Big\}, u_t^i-u_t^{i,o} \ra    \geq  0,  \nonumber \\
&  \forall u^i \in {\mathbb A}^i,  a.e. t \in [0,T], {\mathbb P}|_{{\cal G}_{0,t}^{I^{i,u^o}}}- a.s., i=1,2,\ldots N.   \label{eqh35o} 
\end{align}

\end{description}
Under the additional Assumptions~\ref{a-nf1r} the results also hold for strategies ${\mathbb U}^{(N),y^u}[0,T]$ with conditional expectation taken with respect to ${\cal G}_{0,t}^{y^{i,u^o}}$. \\
For distributed filtering strategies are $u \in {\mathbb U}^{(N),y}[0,T]$ with conditional expectation taken with respect to ${\cal G}_{0,t}^{y^{i,u^o}}$.

\end{theorem}

 \begin{proof}  The derivation consists of two steps. The first step utilizes \cite{charalambous-ahmedFIS_Parti2012} to derive the optimality conditions  for the augmented system (\ref{pds20}). Hence, by direct application of \cite{charalambous-ahmedFIS_Parti2012} we have the following.\\
  Define the Hamiltonian of the augmented system (\ref{pds20}) 
\bes
 {\cal  H}: [0, T]  \times  {\mathbb R}^{n+k} \times {\mathbb R}^{n+k}   \times {\cal L}( {\mathbb R}^{n+k},   {\mathbb R}^{n+k}) \times {\mathbb A}^{(N)}  \longrightarrow {\mathbb R}
\ees  
   by  
   \begin{align}
    {\cal H} (t,X,\Psi,M,u) \tri    \langle F(t,X,u),\zeta \rangle + tr (M^* G(t,X))
     + L(t,X,u),  \hst  t \in  [0, T]. \label{h1piaaa}
    \end{align}
 \noi Then the result of \cite{charalambous-ahmedFIS_Parti2012} for nonrandomized strategies    apply to the system (\ref{pds20}),  hence for any $u \in \times_{i=1}^N {\mathbb U}^{y^{i,u}}[0,T]$ or $u \in \times_{i=1}^N {\mathbb U}^{I^{i,u}}[0,T]$    the adjoint process of the augmented system exists and satisfies the following backward stochastic differential equation.
\begin{align} 
d\Psi (t)  &= -F_X^{*}(t,X(t),u_t)\Psi (t)  dt - V_{Q}(t) dt -L_X(t,X(t),u_t) dt + Q(t) \left[ \begin{array}{c} dW(t) \\ dB(t) \end{array} \right], \hso t \in [0,T),    \nonumber    \\ 
&=- {\cal H}_X (t,X(t),\Psi(t),Q(t),u_t) dt + Q(t) \left[ \begin{array}{c} dW(t) \\ dB(t) \end{array} \right],  \hso \Psi(T)=  \Phi_X(X(T)),  t \in [0,T), \label{h2pi}  
 \end{align}
  where  $V_{Q}$ is   given by  $\langle V_{Q}(t),\zeta\rangle = tr (Q^*(t) G_X(t,X(t); \zeta)), t \in [0,T]$. \\
 The second step translates the necessary conditions of the augmented system  to the original system (\ref{eq1}), (\ref{eq1b}). To this end, we introduce the following decompositions which will lead to a simplified Hamiltonian system of equations.
 \bea
 \Psi \tri \left[ \begin{array}{c} \psi \\ \zeta \end{array} \right], \hst Q \tri \left[ \begin{array}{cc} q_{11} & q_{12} \\ q_{21} & q_{22} \end{array} \right]. \label{dec1}
 \eea
By utilizing this decomposition it can be shown that  $\psi$ satisfies (\ref{adj1a}), (\ref{eq18}).  The second component of $\Psi$ in (\ref{dec1})  satisfies the following equation  
\bea
  d \zeta(t) &=q_{21}(t) dW(t)+ q_{22}(t) dB(t), \hst \zeta(T)=0, \hst t \in [0, T). \label{ph2pia}
    \eea
Since this equation has terminal condition  $\zeta(T)=0$, and its right hand side martingale terms  are orthogonal, then necessarily, $q_{21}(t)=0, q_{22}(t)=0, \forall t \in [0,T], a.s.$, which imply  $\zeta(t)=0, \forall t \in [0,T], a.s$. Finally,  statements {\bf (1)-(3)}  are obtained  from equivalent statements of the augmented system \cite{charalambous-ahmedFIS_Parti2012}.
\end{proof}

 It is interesting to note that the necessary conditions for a $u^o \in {\mathbb U}^{(N),y^u}[0,T]$ or  $u^o \in {\mathbb U}^{(N),I^u}[0,T]$  to be  a person-by-person optimal  can be derived following the methodology of Theorem~\ref{theorem5.1o}, and that these necessary conditions are the same as the necessary conditions for the team optimal strategy. 
These results are  stated as a Corollary.

 \begin{corollary} (Necessary conditions for person-by-person optimality)
 \label{corollary5.1}
   Consider Problem~\ref{problemdet2} under the assumptions of Theorem~\ref{theorem5.1o}.  For  an element $ u^o \in {\mathbb U}^{(N),I^u}[0,T]$ with the corresponding solution $x^o \in B_{{\mathbb F}_T}^{\infty}([0,T], L^2(\Omega,{\mathbb R}^n))$ to be a person-by-person optimal strategy, it is necessary  that 
the statements of Theorem~\ref{theorem5.1o}, {\bf (1), (3)}  hold and statement {\bf (2)} is replaced by
 
   \begin{description}
 \item[(2') ]  The variational inequalities are satisfied:
\begin{align}    
  {\mathbb E} \int_{0}^{T} \la   {\mathbb H}_{u^i} (t,x^o(t), \psi^o(t), q_{11}^o(t), u_t^{o}),u_t^i-u_t^{i,o} \ra dt  
  \geq 0, \hst \forall u^i \in {\mathbb U}^{i,I^{i,u}}[0,T], \hso \forall i \in {\mathbb Z}_N. \label{eq16c}
\end{align}

\end{description}
Under the additional Assumptions~\ref{a-nf1r} the results also hold for strategies ${\mathbb U}^{(N),y^u}[0,T]$ with conditional expectation taken with respect to ${\cal G}_{0,t}^{y^{i,u^o}}$. \\
For distributed filtering strategies are $u \in {\mathbb U}^{(N),y}[0,T]$ with conditional expectation taken with respect to ${\cal G}_{0,t}^{y^{i,u^o}}$.
\end{corollary}

 \begin{proof}  The derivation is based on the procedure of Theorem~\ref{theorem5.1o}, which is completely described   in \cite{charalambous-ahmedFIS_Parti2012} 
\end{proof}

The following remark helps identifying the martingale term in the adjoint process.

\begin{remark}
\label{mart}
According to \cite{charalambous-ahmedFIS_Parti2012},  the Riesz representation theorem for Hilbert space martinagles,  determine the maritingale term of the adjoint process $M_t = \int_{0}^t  \Psi_X^o(s) G(s,X^o(s)) \left[ \begin{array}{c} dW(s)\\ dB(s) \end{array} \right]$, dual to the first martingale term in the variational equation (\ref{eq9}),  hence $Q$ in the adjoint equation, is identified as $ Q(t) \equiv \Psi_X(t) G(t,X(t))$. By translating this to the original system then $q_{11} = \psi_x \sigma, q_{12}=\psi_y D^{\frac{1}{2}}$, provided the derivatives $\psi_x, \psi_y$ exist. 
\end{remark}

\subsection{Sufficient Conditions of Optimality}
\label{sufficient}
In this section, we show that the necessary condition of optimality (\ref{eqh35o}) is also a sufficient condition for optimality, under a convexity conditions on the Hamiltonians and the terminal condition.\\

 \begin{theorem} (Sufficient conditions for team optimality)
 \label{theorem5.1s}
   Consider Problem~\ref{problemdet1}  with strategies from ${\mathbb U}^{(N), I^u}[0,T]$ (respectively  ${\mathbb U}^{(N),y^u}[0,T]$), under Assumptions~\ref{NCD1}, and ${\mathbb A}^i$ a closed, bounded and convex subset of ${\mathbb R}^{k_i}, i=1, \ldots N$.
   Let $(x^o(\cdot), u^o(\cdot))$ denote an admissible state and decision  pair and let $\psi^o(\cdot)$ the corresponding adjoint processes. \\
   Suppose the 
 following conditions hold.
   
\begin{description}

\item[(C4)]  ${\mathbb H} (t, \cdot,\psi,q_{11},\cdot),   t \in  [0, T]$ is convex in $(x,u) \in {\mathbb R}^n \times {\mathbb A}^{(N)}$; 
 
 \item[(C5)] $\varphi(\cdot)$ is convex in $x \in {\mathbb R}^n$.  
\end{description}

\noi Then $(x^o(\cdot),u^o(\cdot))$ is a team optimal pair if it satisfies (\ref{eqh35o}) (respectively (\ref{eqh35o}) with conditional expectation taken in terms of  ${\cal G}_{0,t}^{y^{i,u^o}}, i=1, \ldots N$).
\end{theorem}

 \begin{proof}  Let $u^o \in {\mathbb U}^{(N),y^u}[0,T]$ denote a candidate for the optimal team decision and $u \in {\mathbb U}^{(N),y^u}[0,T]$
any other decision.  Then
 \begin{align} 
 J(u^o) -J(u)=  
    {\mathbb E} \biggl\{   \int_{0}^{T}  \Big(\ell(t,x^o(t),u_t^{o})  -\ell(t,x(t),u_t) \Big)  dt  
     + \Big(\varphi(x^o(T)) - \varphi(x(T))\Big)  \biggr\}    . \label{s1}
  \end{align} 
By the convexity of $\varphi(\cdot)$ then 
\bea
\varphi(x(T))-\varphi(x^o(T)) \geq \la \varphi_x(x^o(T)), x(T)-x^o(T)\ra . \label{s2}
\eea
Substituting (\ref{s2}) into (\ref{s1}) yields
 \begin{align} 
 J(u^o) -J(u)\leq  {\mathbb E} \Big\{ & \la \varphi_x(x^o(T)),   x^o(T) - x(T)\ra \Big\} \nonumber \\
 + &    {\mathbb E} \biggl\{    \int_{0}^{T}  \Big(\ell(t,x^o(t),u_t^{o})  -\ell(t,x(t),u_t) \Big)  dt   \biggr\}    . \label{s3}
  \end{align} 
Applying the Ito differential rule to $\la\psi^o,x-x^o\ra$ on the interval $[0,T]$ and then taking expecation we obtain the following equation.
\begin{align}
{\mathbb E} \Big\{ & \la \psi^o(T),   x(T) - x^o(T)\ra \Big\}
 =  {\mathbb E} \Big\{  \la \psi^o(0),   x(0) - x^o(0)\ra \Big\} \nonumber \\
& +{\mathbb E} \Big\{ \int_{0}^{T}  \la-f_x^{*}(t,x^o(t),u_t^{o})\psi^o(t)dt-V_{q_{11}^o}(t)-\ell_x(t,x^o(t),u_t^{o}), x(t)-x^o(t)\ra dt  \Big\} \nonumber \\
&+ {\mathbb E} \Big\{   \int_{0}^{T}  \la \psi^o(t), f(t,x(t),u_t)- f(t,x^o(t),u_t^{o})\ra dt \Big\} \nonumber \\
& + {\mathbb E} \Big\{  \int_{0}^{T}     tr (q_{11}^{*,o}(t)\sigma(t,x(t))-q_{11}^{*,o}(t)\sigma(t,x^o(t))dt \Big\} \nonumber \\
&=  -  {\mathbb E} \Big\{ \int_{0}^{T} \la {\mathbb H}_x(t,x^o(t),\psi^o(t), q_{11}^o(t),u_t^{o}), x(t)-x^o(t)\ra dt \nonumber \\
& + {\mathbb E} \Big\{  \int_{0}^{T} \la \psi^o(t), f(t,x(t),u_t)-f(t,x^o(t),u_t^{o})\ra dt \Big\} \nonumber \\
& + {\mathbb E} \Big\{  \int_{0}^{T}  tr (q_{11}^{*,o}(t)\sigma(t,x(t))-q_{11}^{*,o}(t)\sigma(t,x^o(t)))dt \Big\} \label{s4}
\end{align}
Note that $\psi^o(T)=\varphi_x(x^o(T))$. Substituting (\ref{s4}) into (\ref{s3}) we obtain
 \begin{align} 
 J(u^o) -J(u)  \leq &    {\mathbb E} \Big\{ \int_{0}^{T}  \Big[ {\mathbb H}(t,x^o(t),\psi^o(t), q_{11}^o(t),u_t^{o}) -   {\mathbb H}(t,x(t),\psi^o(t), q_{11}^o(t),u_t^{})\Big]dt \Big\} \nonumber \\
 -&    {\mathbb E} \Big\{ \int_{0}^{T} \la {\mathbb H}_x(t,x^o(t),\psi^o(t), q_{11}^o(t),u_t^{o}), x^o(t)-x(t)\ra dt   \Big\}    . \label{s5}
  \end{align} 
Since by hypothesis $ {\mathbb H}$ is convex in $(x,u) \in {\mathbb R}^n \times {\mathbb A}^{(N)}$,  then 
\begin{align}
 {\mathbb H} (t,x(t),&\psi^o(t), q_{11}^o(t),u_t^{}) -   {\mathbb H}(t,x^o(t),\psi^o(t), q_{11}^o(t),u_t^{o})  \nonumber \\
  \geq 
 & \sum_{i=1}^N \la {\mathbb H}_{u^i}(t,x^o(t),\psi^o(t), q_{11}^o(t),u_t^{o}), u^i-u_t^{i,o}\ra \nonumber \\
 + & \la {\mathbb H}_x(t,x^o(t),\psi^o(t), Q^o(t),u_t^{o}), x(t)-x^o(t)\ra , \hst t \in [0,T] \label{s7}
 \end{align}
 Substituting (\ref{s7}) into (\ref{s5}) yields 
 \begin{align}
  J(u^o) -J(u)  \leq    
 -    {\mathbb E}  \Big\{  \sum_{i=1}^N    \int_{0}^{T} \la {\mathbb H}_{u^i}(t,x^o(t),\psi^o(t), q_{11}^o(t),u_t^{o}), u_t^i-u_t^{i,o}) dt   \Big\}    . \label{s5a}
  \end{align} 
By     (\ref{eqh35o})  
  and  by definition of conditional expectation   we have
\begin{align}
&{\mathbb E}\Big\{ I_{A_t^i}(\omega)  \la {\mathbb H}_{u^i}(t,x^o(t),\psi^o(t), q_{11}^o(t),u_t^{o}),u_t^i-u_t^{i,o} \ra \Big\}  \nonumber \\
&= {\mathbb E} \Big\{ I_{A_t^i}(\omega) {\mathbb E} \Big\{ \la  {\mathbb H}_{u^i}(t,x^o(t),\psi^o(t),q_{11}^o(t),u_t^{o}),u_t^i-u_t^{i,o} \ra |{\cal G}_{0, t}^{y^{i,u}} \Big\}   \Big\} \geq 0, \hso \forall A_t^i \in {\cal G}_{0,t}^{y^{i,u}},    \hso \forall i \in {\mathbb Z}_N. \label{s5b}
\end{align}
Hence, $\la {\mathbb H}_{u^i}(t,x^o(t),\psi^o(t), q_{11}^o(t),u_t^{o}),u_t^i-u_t^{i,o} \ra\geq 0, \forall u_t^i \in {\mathbb A}^i,  a.e. t \in [0,T], {\mathbb P}- a.s., i=1,2,\ldots, N$. Substituting the this inequality into (\ref{s5a}) gives 
\bes
 J(u^o) \leq J(u), \hst \forall u \in {\mathbb U}^{(N),y^u}[0,T].
\ees
Hence, sufficiency of (\ref{eqh35o})  with conditional expectation taken in terms of  ${\cal G}_{0,t}^{y^{i,u^o}}, i=1, \ldots N$ is shown. For $ {\mathbb U}^{(N),I^u}[0,T]$ the derivation is identical. 
\end{proof} 

Since the necessary conditions for team optimal and person-by-person optimal are equivalent (this follows from Theorem~\ref{theorem5.1o},  Corollary~\ref{corollary5.1}), then one can go one step further to show that under the conditions of Theorem~\ref{theorem5.1s}, that  any person-by-person optimal strategy is also a team optimal strategy.  \\

We conclude our discussion on team and person-by-person game optimality conditions for distributed stochastic differential systems with decentralized noisy information structures, by stating once again that the results derived are also applicable to distributed estimation problems (see Remark~\ref{remf}) with strategies taken from      ${\mathbb U}^{(N),y^u}[0,T]$.

\section{Applications in Communication, Filtering and Control}
\label{classes}
In this section we investigate various applications of the optimality conditions to 
 communication, filtering and control applications. For most applications we give explicit optimal team strategies,  when the  dynamics and the reward have the structures  defined below. Throughout, we assume validity of the convexity conditions of Theorem~\ref{theorem5.1s}, {\bf (C4), (C5)}, why imply sufficiency of (\ref{eqh35o}) with conditional expectation taken in terms of  ${\cal G}_{0,t}^{y^{i,u^o}}, i=1, \ldots N$, and strategies  taken from ${\mathbb U}^{(N),y^u}[0,T]$. \\

\begin{definition}(Team games with special structures)
\label{normal}
We define the following classes of team games.

\textbf{(NF): Nonlinear Form.} The team game is said to have "nonlinear form" if 
\begin{align}
f(t,x,u) \tri &  b(t,x)+ g(t,x)u, \hst g(t,x)u \tri \sum_{j=1}^N g^{(j)}(t,x) u^j,  \label{gn1} \\
 \sigma(t,x)  \tri & \left[ \begin{array}{cccc}  \sigma^{(1)}(t,x) & \sigma^{(2)}(t,x)  \ldots & \sigma^{(N)}(t,x)  \end{array} \right]  \label{gn2} \\
 \ell(t,x,u) \tri &\frac{1}{2}\la u,R(t,x)u\ra + \frac{1}{2} \lambda(t,x)+ \la u, \eta(t,x)\ra , \label{gn3} 
\end{align}
where
\bea
 \la u, R(t,x) u \ra \tri  \sum_{i=1}^N \sum_{j=1}^N u^{i,*} R_{ij}(t,x)u^j, \hst \la u,\eta(t,x)\ra \tri  \sum_{i=1}^N u^{i,*}\eta^i(t,x), \nonumber
\eea
and   $\sigma^{(i)}(\cdot,\cdot)$ is the $i$th column of an $n\times m$ matrix $\sigma(\cdot,\cdot)$, for $i=1,2, \ldots, m$, $R(\cdot,\cdot)$ is symmetric uniformly positive definite, and $\lambda(\cdot,\cdot)$ is uniformly positive semidefinite.

\textbf{(LQF): Linear-Quadratic  Form.} A team game is  said  to have "linear-quadratic form" if     
\begin{align}
f(t,x,u)=&A(t)x+  B(t)u,  \hst \sigma(t,x,u)= G(t), \label{n1} \\
 \ell(t,x)=&\frac{1}{2}\la u,R(t)u\ra  + \frac{1}{2} \la x,H(t)x\ra +\la x,F(t)\ra + \la u, E(t)x\ra  +\la u,m(t)\ra , \label{n3}
\end{align}
and $R(\cdot)$ is symmetric uniformly positive definite and $H(\cdot)$ is symmetric uniformly positive semidefinite.

\end{definition}

\noi Below we compute the optimal strategies for the two cases of Definition~\ref{normal}. First, we introduce the following definitions.
\bes
\widehat{u^{i,j,o}}(t) \tri {\mathbb E}\Big( u_t^{i,o} | {\cal G}_{0,t}^{y^{j,u^o}}\Big), \hso  \widehat{u^{i,o}}(t) \tri Vector\{ \widehat{u^{1,i,o}}(t), \ldots,  \widehat{u^{N,i,o}}(t)\}, i,j=1 \ldots, N,
\ees
\bes
\widehat{{ u}^o}(t) \tri  Vector\{ \widehat{u^{1,o}}(t), \ldots,  \widehat{u^{N,o}}(t)\}, \hso \widehat{x^o}(t) \tri  Vector\{ {\mathbb E}\Big(x^{o}(t) | {\cal G}_{0,t}^{y^{1,u^o}}\Big), \ldots,  {\mathbb E}\Big(x^{o}(t) | {\cal G}_{0,t}^{y^{N,u^o}}\Big)  \},
\ees
\bes
\widehat{\psi^o}(t) \tri  Vector\{ {\mathbb E}\Big(\psi^{o}(t) | {\cal G}_{0,t}^{y^{1,u^o}}\Big), \ldots,  {\mathbb E}\Big(\psi^{o}(t) | {\cal G}_{0,t}^{y^{N,u^o}}\Big) \} , 
\ees
\bes
R^{[i]}(t)\tri \left[ \begin{array}{c} R_{i1}(t), \ldots, R_{iN}(t) \end{array} \right], \hst E^{[i]}(t) \tri \left[ \begin{array}{c} E_{i1}(t), \ldots, E_{iN}(t) \end{array} \right], \hso i=1, \ldots, N.
\ees
\ \

\textbf{Case NF}. \\
\noi Utilizing the definition of Hamiltonian of Theorem~\ref{theorem5.1o},  its derivative is given by
\bea
 {\cal H}_u(t,x,\psi,q_{11},u)  = g^*(t,x) \psi   +  R(t,x)u + \eta(t,x), \hso   (t,x)\in [0,T]\times {\mathbb R}^n. \label{dh3}
 \eea
 The explicit expression for $u_t^{i,o}$ is given by
 \begin{align}
u_t^{i,o}=&- \Big\{  {\mathbb E} \Big( R_{ii}(t,x^o(t))  | {\cal G}_{0,t}^{y^{i,u^o}}\Big)\Big\}^{-1} \Big\{ {\mathbb E} \Big( \eta^i(t,x^o(t))  | {\cal G}_{0,t}^{y^{i,u^o}}\Big)  
 + \sum_{j=1,j \neq i}^N  {\mathbb E} \Big( R_{ij}(t,x^o(t))u_t^{j,o}     | {\cal G}_{0,t}^{y^{i,u^o}} \Big)    \nonumber     \\
  & -     {\mathbb  E} \Big(  g^{(i),*}(t,x^o(t))\psi^{o}(t) |{\cal G}_{0,t}^{y^{i,u^o}}\Big) \Big\}, \hst  {\mathbb P}|_{ {\cal G}_{0,t}^{y^{i,u^o}}}-a.s., \hso i=1,2, \ldots, N .\label{sgrn4}
\end{align}
{\it Special Case.} Suppose $R(t,x)= \overline{R}(t)$, e.g., independent of $x$. Since both sides of (\ref{sgrn4}) are $  {\cal G}_{0,t}^{y^{i,u^o}}-$measurable taking conditional expectations of both side with respect to $  {\cal G}_{0,t}^{y^{i,u^o}}-$ gives the expression
\begin{align}
\widehat{u^{i,i,o}}(t)=&- \Big\{  \overline{R}_{ii}(t)  \Big\}^{-1} \Big\{ {\mathbb E} \Big( \eta^i(t,x^o(t))  | {\cal G}_{0,t}^{y^{i,u^o}}\Big)  
 + \sum_{j=1,j \neq i}^N \overline{R}_{ij}(t) \widehat{u^{j,i,o}}(t)     \nonumber     \\
  & -     {\mathbb  E} \Big(  g^{(i),*}(t,x^o(t))\psi^{o}(t) |{\cal G}_{0,t}^{y^{i,u^o}}\Big) \Big\}, \hst  {\mathbb P}|_{ {\cal G}_{0,t}^{y^{i,u^o}}}-a.s., \hso i=1,2, \ldots, N .\label{sgrn4a}
\end{align}
The last equation can be written in the form of a fixed point matrix equation with random coefficients. We discuss this below. \\

\textbf {Case LQF.} \\
For a  team game is of normal form then from the previous optimal strategies one obtains
\begin{align}
u_t^{i,o}=&-  R_{ii}^{-1}(t) \Big\{  m^i(t)+ \sum_{j=1}^NE_{ij}(t) {\mathbb E} \Big( x^{j,o}(t) | {\cal G}_{0,t}^{y^{i,u^o}}\Big) + \sum_{j=1,j \neq i}^N  R_{ij}(t)  {\mathbb E} \Big( u_t^{j,o} | {\cal G}_{0,t}^{y^{i,u^o}} \Big)    \nonumber     \\
  & -   B^{(i),*}(t)      {\mathbb  E} \Big( \psi^{o}(t) |{\cal G}_{0,t}^{y^{i,u^o}}\Big) \Big\}, \hst  {\mathbb P}|_{ {\cal G}_{0,t}^{y^{i,u^o}}}-a.s., \hso i=1,2, \ldots, N .\label{n4r}
\end{align}
Similarly as above,  (\ref{n4r}) can be put in the form of fixed point  matrix equation as follows. 
\begin{align}
 \diag\{R^{[1]}(t), \ldots, R^{[N]}(t)\} \widehat{{ u}^o}(t) +& \diag\{E^{[1]}(t), \ldots, E^{[N]}(t)\} \widehat{{x}^o}(t) \nonumber \\
 &+ \diag\{B^{(1),*}(t), \ldots, B^{(N),*}(t) \} \widehat{{\bf \psi}^o}(t)+ m(t) =0. \label{fp1}
 \end{align}
Therefore,  (\ref{fp1}) can be solved via  fixed point methods.  One can proceed further to determine the adjoint processes and the explicit optimal team strategy.  This is done in the next subsection.

\subsection{Communication Channels with Memory and Feedback}
\label{static-feedback}
In this section we discuss applications of team games to communication channels with feedback and memory. We consider applications in which the state process is a  RV, and decentralized information structures with feedback and/or  correlation  among them. Consider a filtered probability space $\Big(\Omega, {\mathbb F}, {\mathbb F}_T, {\mathbb P}\Big)$ on which the following are defined.

\begin{align}
&\mbox{A Gaussian RV $ \theta \tri Vector\{\theta^1,\ldots, \theta^N\}: \Omega \longrightarrow {\mathbb R}^n, \theta^i \in {\mathbb R}^{n_i}, ({\mathbb E}(\theta), Cov(\theta))=(\bar{\theta}, P_0)$}, \nonumber  \\
&\mbox{Mutual Independent Brownian motions}  \hso B^i: [0,T] \times \Omega \rar {\mathbb R}^{k_i}, \hso i=1,\ldots N, \: \mbox{independent of $\theta$}. \nonumber 
\end{align}
The information structure of each DM $u^i$ is  ${\cal G}_{0,t}^{y^i} \tri \sigma \{y^i(s): 0\leq s \leq t\}, t \in [0,T]$, $i=1,\ldots, N$, which is defined by a communication channel with memory feedback via the  stochastic differential equation 
\bea
y^i(t) = \int_{o}^t C_{ii}(s,y^i(s))\theta ds +\int_{0}^t D_{ii}^{\frac{1}{2}}(s)dB^i(s), \hst t \in [0,T], \hst i=1,2, \ldots, N, \label{d51a}
\eea
where $C_{ii}: [0,T] \times  {\mathbb R}^{k_i} \longrightarrow {\cal L}({\mathbb R}^n, {\mathbb R}^{k_i})$. The communication channel (\ref{d51a}) models a Gaussian Broadcast channel in which there is a single transmitter and multiple receivers, $i=1, \ldots, N$. The transmitter wishes to send linear combinations of messages $\{\theta^1, \ldots, \theta^N\}$ to receivers $\{y^1, \ldots, y^N\}$. The objective is to reconstruct at each receiver $y^i$ the intended linear combination of the messages denoted by $L^i \theta$, where $L^i$ is an appropriately chosen matrix. A reasonable pay-off for reconstructing the intended linear combination $L^i \theta$ at receiver $i$ by $u^i$ is the average weighted estimation error   ${\mathbb E} \int_{[0,T]}\la u_t- \diag\{L^1, \ldots, L^N\}\theta, R(t)(u_t-   \diag\{L^1, \ldots, L^N\}  \theta)\ra dt$. A more general pay-off which also incorporates any power constraints at the transmitter  is the quadratic pay-off   defined by
\begin{align}
J(u^1, \ldots, u^N)\tri & \frac{1}{2} {\mathbb E} \Big\{ \int_{0}^T \Big( \la u_t,R(t)u\ra  +  \la \theta,H(t)\theta\ra +\la \theta,F(t)\ra  \nonumber \\
&+ \la u_t, E(t) \theta \ra  +\la u_t,m(t)\ra \Big)dt  . \label{d52}
\end{align} 
Noticed that the information structures (\ref{d51a}) are defined via channels with feedback since 
\bea
{\mathbb P} \Big\{ y^i(t) \in A_i | \{y^i(s) : 0 \leq s \leq t-\epsilon\}, \theta \Big\} \neq {\mathbb P} \Big\{ y^i(t) \in A_i | \theta\Big\}, \: A_i \in {\cal B}({\mathbb R}^{k_i}), \eps>0,  i \in {\mathbb Z}_N.
\eea
  The previous communication model can be easily generalized to other network communication channels. \\
Since minimizing (\ref{d52}) over feedback information structures subject to (\ref{d51a}) is a team problem,  then we will apply the optimality conditions of Theorem~\ref{theorem5.1o}. \\
  First, note that   the stochastic differential equation (\ref{d51a}) has a continuous strong solution which is unique. Since the state is a RV (static state), then $\psi^o=0$,  hence  the optimal strategies are given component wise by
\begin{align}
u_t^{i,o}=-  R_{ii}^{-1}(t) \Big\{  \sum_{j=1}^N E_{ij}(t) {\mathbb E} \Big( \theta^j | {\cal G}_{0,t}^{y^{i}}\Big) +m^i(t) + \sum_{j=1,j \neq i}^N  R_{ij}(t)  {\mathbb E} \Big( u_t^{j,o} | {\cal G}_{0,t}^{y^{i}} \Big)  \Big\}, \hst  {\mathbb P}|_{ {\cal G}_{0,t}^{y^{i}}}-a.s., i \in {\mathbb Z}_N .\label{d53}
\end{align}
Define the filter version of $\theta$ by $\widehat{\theta}^i(t) \tri {\mathbb E}\Big( \theta |{\cal G}_{0,t}^{y^{i,o}}\Big), t \in [0,T], i \in {\mathbb Z}_N$. Then these bank of filters  satisfy the following  stochastic differential equations  
\begin{align}
&d\widehat{\theta}^i(t) = P^i(t,y^i) C_{ii}^{*}(t,y^i(t)) D_{ii}^{-1}(t) \Big(dy^i(t)- C_{ii}(t,y^i(t)) \widehat{\theta}^i(t)dt\Big), \hso \widehat{\theta}^i(0)= \overline{\theta}, \hso t \in (0,T], \hso i \in {\mathbb Z}_N \label{d54} \\
&\dot{P}^i(t,y^i) =-P^i(t,y^i) C_{ii}^{*}(t,y^i(t)) D_{ii}^{-1}(t)C_{ii}(t,y^i(t)) P^i(t,y^i), \hso P(0)=\overline{P}_0, \hst t \in (0,T], \hso i \in {\mathbb Z}_N. \label{d56} 
\end{align}
Define the innovations process and the $\sigma-$algebra generated by it as follows.
\begin{align}
I^i(t) & \tri   \int_{0}^t D_{ii}^{\frac{1}{2},-1}(s)\Big( y^i(s)- C_{ii}(s,y^i(s)) \widehat{\theta}^i(s)ds\Big),  \hst  t \in [0,T], \hso i \in {\mathbb Z}_N, \label{d57} \\
  {\cal G}_{0,t}^{I^i} & \tri \sigma\Big\{I^i(s): 0\leq s \leq t\Big\}, \hst  t \in [0,T], \hso i \in {\mathbb Z}_N. \label{d57b}
\end{align}
Then $\{I^i(t): 0\leq t \leq T\}$ is an $\Big({\cal G}_T^{y^i}, {\mathbb P}\Big)-$adapted Brownian motion $\forall i \in {\mathbb Z}_N$, and for $i \neq j$, the innovations $I^i(\cdot), I^j(\cdot)$ are independent (in view of independence of $B^i(\cdot), B^j(\cdot)$ for $i\neq j$).  Moreover, the processes $\{\widehat{\theta}^i, P^i(t,y^i), y^i(t): 0\leq t \leq T\}$ are weak solutions \cite{liptser-shiryayev1977} of the system 
\begin{align}
&d\widehat{\theta}^i(t) = P^i(t,y^i) C_{ii}^{*}(t,y^i(t)) N_{ii}^{-1}(t) dI^i(t), \hst \widehat{\theta}^i(0)= \overline{\theta}, \hst t \in (0,T], \hso i \in {\mathbb Z}_N \label{d58} \\
&dy^i(t) =  C_{ii}(t,y^i(t))\widehat{\theta}^i(t)dt + D_{ii}^{\frac{1}{2}}(t)dI^i(t), \hst t \in [0,T], \hst i \in {\mathbb Z}_N, \label{d51} \\
&\dot{P}^i(t,y^i) =-P^i(t,y^i) C_{ii}^{*}(t,y^i(t)) D_{ii}^{-1}(t)C_{ii}(t,y^i(t)) P^i(t,y^i), \hst P(0)=\overline{P}_0, \hst t \in (0,T], \hso i \in {\mathbb Z}_N. \label{d60} 
\end{align}
Next, we establish existence of strong solutions to the system (\ref{d58})-(\ref{d60}) which will imply that ${\cal G}_T^{y^i}$ and ${\cal G}_T^{I^i}$, $\forall i \in {\mathbb Z}_N$ generate the same information. Under assumption that $C_{ii}(t,y^i)$ satisfy the Lipschitz and linear growth conditions, uniformly in $t \in [0,T]$, the system (\ref{d58}), (\ref{d51}) has a unique  ${\cal G}_T^{I^i}-$adapted continuous solution \cite{liptser-shiryayev1977}, hence $\widehat{\theta}^i(\cdot)$     is  ${\cal G}_{0,t}^{I^i}-$measurable, $\forall t \in [0,T]$,  $\forall i \in {\mathbb Z}_N$. Thus, ${\cal G}_{0,t}^{y^i} \subseteq {\cal G}_{0,t}^{I^i}, \forall t \in [0,T]$, $\forall i \in {\mathbb Z}_N$. The reverse ${\cal G}_{0,t}^{I^i} \subseteq {\cal G}_{0,T}^{y^i}, \forall t \in [0,T] $,  $\forall i \in {\mathbb Z}_N$ follows from the construction of innovations processes (\ref{d57}). Hence, ${\cal G}_T^{I^i} = {\cal G}_T^{y^i} $,  $\forall i \in {\mathbb Z}_N$. Since each DM $u^i$ is ${\cal G}_T^{y^i}={\cal G}_T^{I^i}-$adapted $\forall i \in {\mathbb Z}_N$, and the innovations sigma algebras ${\cal G}_T^{I^i}$ are independent for $\forall i, j \in {\mathbb Z}_N, i\neq j$ then the optimal strategies (\ref{d53}) are given  by
\begin{align}
u_t^{i,o}=-  R_{ii}^{-1}(t) \Big\{   E^{[i]}(t) \widehat{\theta}^i(t)+m^i(t) + \sum_{j=1,j \neq i}^N  R_{ij}(t)  {\mathbb E} \Big( u_t^{j,o}  \Big)  \Big\}, \hst  {\mathbb P}|_{ {\cal G}_{0,t}^{y^i}}-a.s., i \in {\mathbb Z}_N .\label{d61}
\end{align}
Next, we determine the vector by   $\overline{u^o}  \tri Vector\{   
{\mathbb  E}(u^{1,o}), {\mathbb E} (u_t^{2,o}),\ldots,
 {\mathbb E}(u_t^{N,o})\}.$
Taking expectation of both sides of (\ref{d61}) gives the following linear system of equations.
\begin{align}
\overline{u_t^{i,o}}(t)=-  R_{ii}^{-1}(t) \Big\{   E^{[i]}(t) {\mathbb E} \Big( \theta \Big) +m^i(t) + \sum_{j=1,j \neq i}^N  R_{ij}(t)  \overline{u_t^{j,o}}(t)  \Big\}, \hst  i \in {\mathbb Z}_N .\label{d62}
\end{align}
The last equation can be put into a fixed point form.
Define 
\begin{align}
      M(t) \tri \left[ \begin{array}{c}  
-R_{11}^{-1}(t) E^{[1]}(t) {\mathbb E} \Big(\theta\Big) \\
-R_{22}^{-1}(t) E^{[2]}(t) {\mathbb E} \Big(\theta\Big) \\
\ldots \\
\ldots\\
-R_{NN}^{-1}(t) E^{[N]}(t) {\mathbb E} \Big(\theta\Big) 
\end{array} \right],  \hst  K(t) \tri \left[ \begin{array}{c}  
-R_{11}^{-1}(t) m^1(t)  \\
-R_{22}^{-1}(t) m^2(t)  \\
\ldots \\
\ldots\\
-R_{NN}^{-1}(t) m^{N}(t)  
\end{array} \right],     
   \label{d64}
\end{align}
\begin{align}
\Lambda(t) \tri \left[ \begin{array}{ccccc}     I & R_{11}^{-1}(t)R_{12}(t) & R_{11}^{-1}(t) R_{13}(t) & \ldots & R_{11}^{-1}(t) R_{1N}(t) \\ \\
  R_{22}^{-1}(t)R_{21}(t) & I & R_{22}^{-1}(t) R_{23}(t) & \ldots & R_{22}^{-1}(t) R_{2N}(t) \\ \\
\ldots &  \ldots   & I  & \ldots &  \ldots \\ \\
\ldots & \ldots &\ldots & \ldots  & \ldots \\ \\
 R_{NN}^{-1}(t)R_{N1}(t) & R_{NN}^{-1}(t) R_{N2}(t)  & R_{NN}^{-1}(t) R_{N3}(t) & \ldots & I \end{array} \right] \label{d65}
 \end{align}
From (\ref{d62}), we have 
\bea
\Lambda(t) \overline{u^o}(t) =  M(t)+ K(t), \hso
\overline{u^o}(t) = \Lambda^{-1}(t)\Big( M(t)+ K(t)\Big), \: t \in [0,T], \: \mbox{if} \: \Lambda(t) >0, \forall t \in [0,T]. \label{d66} 
\eea
Finally, the optimal strategies are given by  (\ref{d61}) and (\ref{d66}).  

The previous calculations can be  generalized to other channel models. Moreover, $\theta$ can be extended to a Random process described by It\^o stochastic differential equations. For linear dynamics this generalization is a straight forward repetition  of the previous calculations, hence it is omitted.

\subsection{Linear-Quadratic Form and Linear Stochastic Differential Dynamics}
\label{lqg}
In this section we invoke the minimum principle to compute the optimal strategies, with respect to a quadratic pay-off, for  distributed stochastic dynamical decision systems consisting of  an interconnection of two subsystems,  each governed by a linear stochastic differential equation with coupling.

\begin{align}
\mbox{\bf Subsystem Dynamics 1:}&  \nonumber \\
  dx^1(t)=& A_{11}(t)x^1(t)dt + B_{11}(t)u_t^1dt  + G_{11}(t) dW^1(t)   \nonumber \\
  &+ A_{12}(t)x^2(t)dt  + B_{12}(t)u_t^2 dt , \hst   x^1(0)=x^1_0, \hso   t \in (0,T],   \label{ex30} \\ \nonumber \\
\mbox{\bf Subsystem Dynamics 2:}&  \nonumber  \\
dx^2(t)=&  A_{22}(t)x^2(t)dt  + B_{22}(t)u_t^2 dt + G_{22}(t) dW^2(t)  \nonumber \\
&+ A_{21}(t)x^1(t)dt + B_{21}u_t^1dt , \hst x^2(0)=x^2_0, \hso  t \in (0,T] \label{ex31} 
\end{align}
For any $t \in [0,T]$ the feedback  information structure       of  $u_t^1$ of subsystem $1$ is  the $\sigma-$algebra ${\cal G}_{0,t}^{y^{1,u}} \tri \sigma\{y^1(s): 0 \leq s \leq t\}$, and 
 the feedback  information structure of  $u_t^2$ of subsystem $2$ is the $\sigma-$algebra ${\cal G}_{0,t}^{y^{2,u}} \tri \sigma\{y^2(s): 0 \leq s \leq t\}$. These  information structures are  defined  by the following linear observation equations.\\
 
{\bf Information structure of Local Control $u^1$:}\\
\bea
y^1(t) =\int_{0}^{t} C_{11}(s) x^1(s) ds  + \int_{0}^t D_{11}^{\frac{1}{2}}(s) d B^1(s), \hst t \in [0,T]. \label{ex32}
\eea

\  \

{\bf Information structure of Local Control $u^2$:} \\
\bea
y^2(t) = \int_{0}^{t} C_{22}(s) x^2(s) ds + \int_{0}^t D_{22}^{\frac{1}{2}}(s) d B^2(s), \hst t \in [0,T]. \label{ex33}
\eea
We may also assume the DMs strategies $u^1$ and $u^2$ are functionals of the innovations information structures ${\cal G}_{0,t}^{I^{1,u}} \tri \sigma\{I^1(s): 0 \leq s \leq t\}$,    ${\cal G}_{0,t}^{I^{2,u}} \tri \sigma\{I^1(s): 0 \leq s \leq t\}$ defined by the innovations processes of (\ref{ex32}), (\ref{ex33}), respectively.\\

The pay-off or reward is quadratic in $(x^1,x^2,u^1,u^2)$. \\

{\bf Pay-off  Functional: }\\

\begin{align}
J(u^1,u^2) =& \frac{1}{2}  {\mathbb E} \Big\{ \int_{0}^T \Big[ \la \left( \begin{array}{c} x^{1}(t) \\ x^{2}(t) \end{array} \right), H(t)   \left(\begin{array}{c} x^1(t) \\ x^2(t) \end{array} \right)      \ra
+ \la \left( \begin{array}{c} u_t^{1} \\ u_t^{2} \end{array} \right), R(t)   \left(\begin{array}{c} u_t^1 \\ u_t^2 \end{array} \right) \ra         \Big]dt \nonumber \\
& + \la   \left( \begin{array}{c} x^{1}(T) \\ x^{2}(T) \end{array} \right), M(T)   \left(\begin{array}{c} x^1(T) \\ x^2(T) \end{array} \right) \ra            \Big\}.  \label{ex34}
\end{align}
We assume that the initial condition $x(0)$, the system Brownian motion $\{W(t): t \in [0,T]\}$, and the observations Brownian motion $\{B^1(t): t \in [0,T]\}$, and $\{B^2(t): t \in [0,T]\}$ are mutually independent and $x(0)$ is Gaussian $({\mathbb E}(x(0)), Cov(x(0)))=(\bar{x}_0, P_0).$\\
For decentralized filtering we set $u^1=0, u^2=0$ in the right hand sides of (\ref{ex30}), (\ref{ex31}), but we should take as   pay-off (\ref{d52}).\\

Define the augmented variables by

\begin{align}
x \tri \left(\begin{array}{c} x^1 \\ x^2 \end{array} \right), \hso y \tri & \left(\begin{array}{c} y^1 \\ y^2 \end{array} \right),\hso u \tri \left(\begin{array}{c} u^1 \\ u^2 \end{array} \right),  \hso \psi \tri \left(\begin{array}{c} \psi^1 \\ \psi^2 \end{array} \right), \hso q_{11} \tri \left(\begin{array}{c} q_{11}^1 \\ q_{11}^2 \end{array} \right), \nonumber \\
W \tri & \left(\begin{array}{c} W^1 \\ W^2 \end{array} \right),\hso B \tri \left(\begin{array}{c} B^1 \\ B^2 \end{array} \right), \nonumber 
\end{align}
and matrices by
\bes
& & A \tri \left[\begin{array}{cc} A_{11}  & A_{12} \\ A_{21} & A_{22} \end{array} \right],  \hso  B \tri \left[\begin{array}{cc} B_{11}  & B_{12} \\ B_{21} & B_{22} \end{array} \right], \hso C \tri \left[\begin{array}{cc} C_{11}  & 0 \\ 0 & C_{22} \end{array} \right], \\
& &  B^{(1)} \tri \left[ \begin{array}{c} B_{11}  \\ B_{21} \end{array} \right], \hso  B^{(2)} \tri \left[ \begin{array}{c} B_{12} \\ B_{22} \end{array} \right], \hso
 C^{[1]} \tri \left[ \begin{array}{cc} C_{11} & 0 \end{array} \right], \hso  C^{[2]} \tri \left[ \begin{array}{cc} 0 & C_{22} \end{array} \right], \\
& & G\tri \left[\begin{array}{cc} G_{11}  & 0 \\ 0 & G_{22} \end{array} \right],  \hso
  D^{\frac{1}{2}}\tri \left[\begin{array}{cc} D_{11}^{\frac{1}{2}}  & 0 \\ 0 & D_{22}^{\frac{1}{2}} \end{array} \right] .
\ees
\noi The distributed system is described in compact form by 
\begin{align}
dx(t)
=&A(t)x(t)dt + B(t)u_t dt + G(t)dW(t), \hst x(0)=x_0 \hso t \in [0, T], \label{ex1} \\
dy(t) =&C(t) x(t) dt  + D^{\frac{1}{2}}(t) d B(t), \hst t \in [0,T]. \label{exo} 
\end{align}
while the pay-off is expressed by
\bea
J(u^1, u^2) =\frac{1}{2}  {\mathbb E}  \Big\{\int_{0}^T  \Big[\la x(t), H(t)x(t)\ra  + \la u_t, R(t) u_t\ra \Big]dt  +  \la x(T),M(T) x(T)\ra \Big\}. \label{ex2} 
\eea
\noi  By Theorem~\ref{theorem5.1o}, the Hamiltonian is given by
\bea
{\mathbb H}(t,x,\psi, q_{11}, u) = \la A(t)x+Bu, \psi\ra  +tr (q_{11}^* G(t)) + \frac{1}{2} \la x, H(t)x\ra  + \frac{1}{2} \la u, R(t) u\ra . \label{ex3}
\eea
The derivative of the Hamiltonian with respect to $u$ component wise this is given by
\begin{align}
{\mathbb H}_{u^1}(t,x,\psi, Q, u^1,u^2) &= B^{(1),*}(t) \psi(t)+  R_{11}(t)u^1 + R_{12}(t)u^2, \label{ex3ab} \\
{\mathbb H}_{u^2}(t,x,\psi, Q, u^1,u^2) &= B^{(2),*}(t) \psi(t)+  R_{22}(t)u^2 + R_{21}(t)u^1. \label{ex3ac} 
\end{align}
The optimal decision $\{u_t^o=(u_t^{1,o}, u_t^{2.o}): 0\leq t \leq T\}$ is obtained from (\ref{ex3ab}), (\ref{ex3ac}) by using  the information structure available to each DM. 

Let $(x^o(\cdot), \psi^o(\cdot), q_{11}^o(\cdot), q_{12}^o(\cdot))$ denote the solutions of the Hamiltonian system, corresponding to the optimal control $u^o$, then 
\begin{align}
 dx^o(t) =& A(t)x^o(t)dt + B(t)u_t^o dt + G(t)dW(t), \hst x^o(0)=x_0 \label{ex5} \\ \nonumber \\
  d\psi^o(t)=& -A^*(t)\psi^o(t)dt   - H(t) x^o(t) dt  -V_{q_{11}^o}(t) dt \nonumber \\
&+ q_{11}^o(t) dW(t)+ q_{12}^o(t) dB(t), \hst \psi^o(T)=M(T) x^o(T) \label{ex6}
\end{align}
Next, we identify the martingale term in (\ref{ex6}). Let $\{\Phi(t,s): 0\leq s \leq t \leq T\}$ denote the transition operator of $A(\cdot)$ and $\Phi^*(\cdot, \cdot)$ that of the adjoint $A^*(\cdot)$ of $A(\cdot)$.  Then $\{ \psi^o(t): t \in [0,T]\}$ is given by 
\begin{align}
\psi^o(t)=& \Phi^*(T,t)M(T) x^o(T) + \int_{t}^T \Phi^*(s,t) \Big\{ H(s) x^o(s) ds + V_{q_{11}^o}(s)ds \nonumber \\
&- q_{11}^o(s) dW(s) -q_{12}^o(s)dB(s) \Big\}. \label{ex9}
\end{align}  
By using the using the identity $\frac{\partial}{\partial s} \Phi^*(t,s) = -A^*(s) \Phi^*(t,s), 0 \leq s \leq t \leq T$ one can verify by differentiation that (\ref{ex9}) is a solution of $(\psi^o(\cdot), q_{11}^o(\cdot),q_{12}^o(\cdot))$ governed by (\ref{ex6}). Since for any control policy, $\{x^o(s): 0\leq t \leq s \leq T\}$ is uniquely determined from (\ref{ex5}) and its current value $x^o(t)$, then (\ref{ex9}) can be expressed via 
\bea
\psi^o(t)=\Sigma(t) x^o(t)+ \beta^o(t) , \hst t \in [0,T], \label{ex15g}
\eea
where $\Sigma(\cdot), \beta^o(\cdot)$ determine the operators to the one expressed via (\ref{ex9}). \\
Next, we determine the operators $(\Sigma(\cdot), \beta^o(\cdot))$. Applying the It\^o differential rule to both sides of (\ref{ex15g}),  and then using (\ref{ex5}), (\ref{ex6}) we obtain
\begin{align}
-A^*(t)\psi^o(t)dt   &- H(t) x^o(t) dt  -V_{q_{11}^o}(t) dt + q_{11}^o(t) dW(t) +q_{12}^o(t) dB(t) \nonumber \\
 =& \dot{\Sigma}(t)x^o(t) dt +  \Sigma(t) \Big\{ A(t)x^o(t)dt + B(t)u_t^o dt + G(t)dW(t)\Big\} + d\beta^o(t). \label{ex16g}
\end{align}
Substituting the claimed relation (\ref{ex15g}) into (\ref{ex16g}) we obtained the identity 
\begin{align}
\Big\{-A^*(t) \Sigma(t)   &-\Sigma(t) A(t) - H(t) -\dot{\Sigma}(t) \Big\} x^o(t) dt  -V_{q_{11}^o}(t) dt + q_{11}^o(t) dW(t)+ q_{12}^o(t)dB(t) \nonumber \\
 =&  A^*(t) \beta^o(t)dt + \Sigma(t)  B(t)u_t^o dt + \Sigma(t) G(t)dW(t) + d\beta^o(t). \label{ex16ga}
\end{align}
 Since $\sigma(t,x)=G(t)$,  then $V_{q_{11}^o}(t)=0,  \forall t \in [0,T]$. 
 By matching the intensity of the martingale terms $\{\cdot\} dW(t)$ in (\ref{ex16ga}), and the rest of the terms  we obtain the following equations.
 \begin{align}
 &V_{q_{11}^o}(t)= 0, \hst \forall t \in [0,T], \label{ex17v} \\
 &q_{11}^o(t)=\Sigma(t) G(t), \hst t \in [0,T], \label{ex17v1} \\
&\dot{\Sigma}(t) +A^*(t) \Sigma(t)   + \Sigma(t) A(t) + H(t)=0,  \hst t \in [0,T), \hst \Sigma(T)=M(T),  \label{exg17} \\
&d \beta^o(t) +A^*(t) \beta^o(t)dt+ \Sigma(t)  B(t)u_t^o dt- q_{12}^o(t) dB(t)=0, \hst t \in [0,T), \hst\beta^o(T)=0. \label{ex18g}
\end{align}
Notice that $q_{12}^o$ is also obtained  by Remark~\ref{mart}, since $q_{12}^o(t)= \psi_x^o(t) G(t)= \Sigma(t)G(t),  \forall t \in [0,T]$.\\

By Theorem~\ref{theorem5.1o}, $\{(u_t^{1,o}, u_t^{2,o}): 0 \leq t\leq T\}$  obtained from (\ref{ex3ab}) and (\ref{ex3ac}), are given by 
\begin{align}
{\mathbb E} & \Big\{ {\cal H}_{u^1} (t,x^{1,o}(t), x^{2,o}(t), \psi^{1,o}(t),\psi^{2,o}(t), q_{11}^{1,o}(t),q_{11}^{2,o}(t), u_t^{1,o}, u_t^{2,0}) |{\cal G}_{0,t}^{y^{1,u^o}}\Big\} =0, \nonumber \\
& \hst a.e.t \in [0,T], \hso {\mathbb P}|_{{\cal G}_{0,t}^{y^{1,u^o}}}-a.s. \label{ex35}
\end{align}
\begin{align}
{\mathbb E} & \Big\{ {\cal H}_{u^2} (t,x^{1,o}(t), x^{2,o}(t), \psi^{1,o}(t),\psi^{2,o}(t), q_{11}^{1,o}(t),q_{11}^{2,o}(t), u_t^{1,o}, u_t^{2,0}) |{\cal G}_{0,t}^{y^{2,u^o}}\Big\} =0, \nonumber \\
& \hst a.e.t \in [0,T], \hso {\mathbb P}|_{{\cal G}_{0,t}^{y^{2,u^o}}}-a.s. \label{ex36}
\end{align}
where $(x^o(\cdot), \psi^o(\cdot), q_{11}^o(\cdot)) \equiv (x^{1,o}(\cdot), x^{2,o}(\cdot), \psi^{1,o}(\cdot), \psi^{2,o}(\cdot),q_{11}^{1,o}(\cdot)), q_{11}^{2,o}(\cdot)) $ are solutions of the Hamiltonian system (\ref{ex5}), (\ref{ex6}) corresponding to $u^o$.  From (\ref{ex35}), (\ref{ex36}) the optimal decisions are 
\bea
u_t^{1,o} = -R_{11}^{-1}(t)   B^{(1),*}(t)  {\mathbb E} \Big\{ \psi^{o}(t) |{\cal G}_{0,t}^{y^{1,u^o}}\Big\} -R_{11}^{-1}(t)R_{12}(t) {\mathbb E} \Big\{u_t^{2,o} | {\cal G}_{0,t}^{y^{1,u^o}}\Big\}  ,  \hst t \in [0,T]. \label{ex41}
\eea
\bea
u_t^{2,o} = -R_{22}^{-1}(t) B^{(2),*}(t)  {\mathbb E} \Big\{ \psi^{o}(t) |{\cal G}_{0,t}^{y^{2,u^o}}\Big\} -R_{22}^{-1}(t)R_{21}(t) {\mathbb E} \Big\{u_t^{1,o} | {\cal G}_{0,t}^{y^{2,u^o}}\Big\}.  \hst t \in [0,T]. \label{ex42}
\eea
Clearly, the previous equations illustrate the coupling between the two subsystems, since $u^{1,o}$ is estimating the optimal decision of the other subsystem  $u^{2,o}$ as well as the adjoint processes $\psi^{o}$ from its own observations, and vice-versa. 

Let $ \phi(\cdot)$ be any   square integrable and ${\mathbb F}_T-$adapted matrix-valued process or scalar-valued processes,  and define its filtered and predictor  versions by   
\bes
\pi^i (\phi)(t) \tri {\mathbb E} \Big\{ \phi(t) | {\cal G}_{0,t}^{y^{i} } \Big\}, \hst  \pi^i (\phi)(s,t) \tri {\mathbb E} \Big\{ \phi(s) | {\cal G}_{0,t}^{y^{i} }  \Big\}, \hst t \in [0,T], \hso s \geq t, \hso i=1,2.
\ees
For any admissible  decision $u$ and corresponding  $(x(\cdot), \psi(\cdot)$ define their  filter versions  with respect to ${\cal G}_{0,t}^{y^i}$ for $i=1, 2$,  by 
\begin{align}
\pi^i(x)(t)  \tri    \left[ \begin{array}{c}  {\mathbb E} \Big\{x^{1}(t)  |  {\cal G}_{0,t}^{y^{i,u}} \Big\}   \\  {\mathbb E}  \Big\{ x^{2}(t) |  {\cal G}_{0,t}^{y^{i,u}} \Big\}  \end{array}     \right], \hso
 \pi^i(\psi)(t) \tri \left[\begin{array}{c}  {\mathbb E} \Big\{   \psi^{1}(t)  |  {\cal G}_{0,t}^{y^{i,u}} \Big\}   \\   {\mathbb E}  \Big\{ \psi^{2}(t)   |    {\cal G}_{0,t}^{y^{i,u}} \Big\}  \end{array}    \right],  \nonumber 
    \hst t \in [0,T], \hso i=1,2,
    \end{align}
\bes
\pi^i(u)(t)  \tri    \left[ \begin{array}{c}  {\mathbb E} \Big\{u_t^{1}  |  {\cal G}_{0,t}^{y^{i,u}} \Big\}   \\  {\mathbb E}  \Big\{ u_t^{2}|  {\cal G}_{0,t}^{y^{i,u}} \Big\}  \end{array}     \right]   \equiv  \left[ \begin{array}{c}  \widehat{u_t^{1,i}}    \\  \widehat{u_t^{2,i}} \end{array}     \right], \hso t \in [0,T] , \hso i=1,2, 
\ees
and their predictor versions by
\bes
\pi^i(x)(s,t)  \tri    \left[ \begin{array}{c}  {\mathbb E} \Big\{x^{1}(s)  |  {\cal G}_{0,t}^{y^{i,u}} \Big\}   \\  {\mathbb E}  \Big\{ x^{2}(s) |  {\cal G}_{0,t}^{y^{i,u}} \Big\}  \end{array}     \right],   \hso
 \pi^i(\psi)(s,t) \tri \left[\begin{array}{c}  {\mathbb E} \Big\{   \psi^{1}(s)  |  {\cal G}_{0,t}^{y^{i,u}} \Big\}   \\   {\mathbb E}  \Big\{ \psi^{2}(s)   |    {\cal G}_{0,t}^{y^{i,u}} \Big\}  \end{array}    \right], 
    \hst t \in [0,T], \hso s \geq t, \hso i=1,2.
\ees
\bes
\pi^i(u)(s,t)  \tri    \left[ \begin{array}{c}  {\mathbb E} \Big\{u_s^{1}  |  {\cal G}_{0,t}^{y^{i,u}} \Big\}   \\  {\mathbb E}  \Big\{ u_s^{2} |  {\cal G}_{0,t}^{y^{i,u}} \Big\}  \end{array}     \right], \hst t \in [0,T], \hso s\geq t, \hso i=1,2, 
\ees
From (\ref{ex41}), (\ref{ex42}) the optimal decisions are 
\begin{align}
u_t^{1,o} \equiv  -R_{11}^{-1}(t) B^{(1),*}(t) \pi^1(\psi^o)(t)   -R_{11}^{-1}(t)R_{12}(t) {\mathbb E} \Big\{u_t^{2,o} | {\cal G}_{0,t}^{y^{1,u^o}}\Big\},   \hst t \in [0,T],    \label{ex41a}
\end{align}
\begin{align}
u_t^{2,o}  \equiv  -R_{22}^{-1}(t) B^{(2),*}(t) \pi^2(\psi^o)(t)-R_{22}^{-1}(t)R_{21}(t) {\mathbb E} \Big\{u_t^{1,o} | {\cal G}_{0,t}^{y^{2,u^o}}\Big\}, \hst t \in [0,T].  \label{ex42a}
\end{align}
The previous  optimal decisions  require  the conditional estimates \\$\{(\pi^1(\psi^o)(t),    \pi^2(\psi^o)(t)): 0\leq t \leq T\}$. These are obtained by  taking conditional expectations of (\ref{ex9}) giving   
\begin{align}
\pi^i(\psi^{o})(t)= \Phi^*(T,t)M(T) \pi^i(x^{o})(T,t) + \int_{t}^T \Phi^*(s,t)  H(s) \pi^i(x^{o})(s,t)ds,   \hst t \in [0,T], \hso i=1,2. \label{ex49}
\end{align}  
For any admissible  decision, the filtered versions of $x(\cdot)$ are given by the following  stochastic differential equations  \cite{liptser-shiryayev1977}.
\begin{align}
&d\pi^1(x)(t) = A(t) \pi^1(x)(t)dt + B^{(1)}(t) u_t^{1} dt + B^{(2)}(t) \pi^1(u^{2})(t) dt 
+ \Big\{ \pi^1(x x^{*})(t)  \nonumber \\
&-\pi^1(x)(t) \pi^1(x^{*})(t) \Big\} C^{[1],*}  D_{11}^{-1}\Big(dy^{1}(t)- C^{[1]}(t) \pi^1(x)(t) dt \Big), \hso \pi^1(x)(0)=\bar{x}_0, \label{ex50}
\end{align}
\begin{align}
&d\pi^2(x)(t) = A(t) \pi^2(x)(t)dt + B^{(2)}(t) u_t^{2} dt + B^{(1)}(t)\pi^2(u^{1})(t) dt  
+ \Big\{ \pi^2(x x^{*})(t) \nonumber \\ 
&-\pi^2(x)(t) \pi^2(x^{*})(t) \Big\}   C^{[2],*}  D_{22}^{-1}(t) \Big(dy^{2}(t)- C^{[2]}(t) \pi^2(x)(t) dt \Big), \hso \pi^2(x)(0)
=\bar{x}_0. \label{ex51}
\end{align}
From the previous filtered versions of $x(\cdot)$ it is clear that  subsystem $1$ estimates the  actions of subsystem $2$ based on its own observations, namely, $\pi^1(u^{2})(\cdot)$ and subsystem $2$ estimates the actions of subsystem $1$ based on its own observations, namely, $\pi^2(u^{1})(\cdot)$. \\
For any admissible decision $(u^1,u^2) \in {\mathbb U}_{reg}^{(2),y^u}[0,T]$ define the innovation processes associated with $\{{\cal G}_{0,t}^{y^{i,u}}: t \in [0,T]\}, i=1,2$ and the  $\sigma-$algebras generated by them as follows 
\bea
I^{i}(t) \tri y^i(t)- \int_{0}^t  C^{[i]}(s)  \pi^i(x)(s)ds,  \hso {\cal G}_{0,t}^{I^{i,u}} \tri \sigma\Big\{ I^i(s): 0\leq s \leq t\Big\}, \hso  t \in [0,T], \hso i=1,2. \label{in1}
\eea 
Let $I^{i,o}(t)$ the innovations processes corresponding to where $(x^o,u^o), i=1,2$. Then by (\ref{in1}), $\{I^i(t): t \in [0,T]\}$ is $\Big( \{ {\cal G}_{0,t}^{y^{i,u}}: t \in [0,T]\}, {\mathbb P}\Big)-$ adapted Brownian motion, $I^i(t)$ has covariance $Cov(I^i(t)) \tri \int_{0}^t D_{ii}(s)ds$, for $i=1,2$ and $\{I^1(t): t \in [0,T]\}$,   $\{I^2(t): t \in [0,T]\}$ are independent.   \\
For any admissible decision $u$ the predicted versions of $x(\cdot)$ are obtained from (\ref{ex50}) and (\ref{ex51}) as follows. Utilizing the identity $\pi^i(x)(s,t)={\mathbb E}  \Big\{  {\mathbb E} \Big\{ x(s)| {\cal G}_{0,s}^{y^{i,u} }\Big\} |  {\cal G}_{0,t}^{y^{i,u}} \Big\}= {\mathbb E} \Big\{ \pi^i(x)(s) | {\cal G}_{0,t}^{y^{i,u}} \Big\}$, for $0 \leq t \leq s \leq T$ then 
\begin{align}
&d\pi^1(x)(s,t) = A(s) \pi^1(x)(s,t)ds + B^{(1)}(s) \pi^1(u^1)(s,t)ds +B^{(2)}(s) \pi^1(u^{2})(s,t)ds, 
 \hst    t < s \leq T, \label{ex52} \\
&\pi^1(x)(t,t)=\pi^1(x)(t), \hst t \in [0,T), \label{ex52a} 
\end{align}
\begin{align}
&d\pi^2(x)(s,t) = A(s) \pi^2(x)(s,t)ds +B^{(2)}(s) \pi^2(u^{1})(s,t) ds+ B^{(1)}(s) \pi^2(u^{1})(s,t) ds, \hst  t < s \leq T, \label{ex53}  \\
&\pi^2(x)(t,t)=\pi^2(x)(t), \hst t \in [0,T).
\label{ex53a}
\end{align}

Since for a given admissible policy and observation paths, $\{\pi^1(x)(s,t): 0\leq t \leq s \leq T\}$ is  determined from (\ref{ex52}) and its current value $\pi^1(x^o)(t,t)=\pi^1(x)(t)$, and $\{\pi^2(x)(s,t): 0\leq t \leq s \leq T\}$ is  determined from (\ref{ex53}) , and its current value $\pi^2(x)(t,t)=\pi^2(x)(t)$, then (\ref{ex49}) can be expressed via
\bea
\pi^i(\psi^o)(t) = K^i(t) \pi^i(x^o)(t) + r^i(t), \hst t\in [0,T], \hso i=1,2. \label{in54}
\eea
 where $K^i(\cdot), r^i(\cdot)$  determines the operators  to the one expressed via (\ref{ex49}), for $i=1,2$. Utilizing (\ref{in54})  into (\ref{ex41a}) and (\ref{ex42a}) then 
\begin{align}
u_t^{1,o} \equiv  -R_{11}^{-1}(t) B^{(1),*}(t) \Big\{K^1(t) \pi^1(x^o)(t) +r^1(t) \Big\}  -R_{11}^{-1}(t)R_{12}(t) \pi^1({u_t^{2,o}})(t),   \hst t \in [0,T],    \label{in55}
\end{align}
\begin{align}
u_t^{2,o}  \equiv  -R_{22}^{-1}(t) B^{(2),*}(t) \Big\{ K^2(t)    \pi^2(x^o)(t) +r^2(t)\Big\} -R_{22}^{-1}(t)R_{21}(t) \pi^2({u^{1,o}})(t), \hst t \in [0,T].  \label{in56}
\end{align}
Let $\{\Psi_{K^i}(t,s): 0\leq s \leq t \leq T\}$ denote the transition operator of $A_{K^i}(t) \tri \Big(A(t)-B^{(i)}(t)R_{ii}^{-1}(t) B^{(i),*}(t) K^i (t) \Big)$, for $i=1,2$. \\
 Next, we determine $K^i(\cdot), r^i(\cdot), i=1,2$. Substituting  the previous equations into (\ref{ex52}), (\ref{ex52a}) and (\ref{ex53}), (\ref{ex53a}) then
  \begin{align}
 \pi^{1}(x^o)(s,t)=& \Psi_{K^1}(s,t) \pi^1(x^o)(t) -\int_{t}^s \Psi_{K^1}(s,\tau) B^{(1)}(\tau)R_{11}^{-1}(\tau)B^{(1),*}(\tau) r^1(\tau)d\tau  \nonumber \\
 &-\int_t^s \Psi_{K^1}(s,\tau) B^{(1)}(\tau)R_{11}^{-1}(\tau) R_{12}(\tau) \pi^1(u^{2,o})(\tau,t)d\tau  \nonumber \\
 &+ \int_{t}^s \Psi_{K^1}(s, \tau) B^{(2)}(\tau) \pi^1(u^{2,o})(\tau,t) d\tau , \hst t \leq s \leq T, \label{pr1}
 \end{align}
  \begin{align}
 \pi^{2}(x^o)(s,t)=& \Psi_{K^2}(s,t) \pi^2(x^o)(t)  -\int_{t}^s \Psi_{K^2}(s,\tau) B^{(2)}(\tau)R_{22}^{-1}(\tau)B^{(2),*}(\tau) r^2(\tau)d\tau  \nonumber \\
 &-\int_t^s \Psi_{K^2}(s,\tau) B^{(2)}(\tau)R_{22}^{-1}(\tau) R_{21}(\tau) \pi^2(u^{1,o})(\tau,t)d\tau  \nonumber \\
 &+ \int_{t}^s \Psi_{K^2}(s, \tau) B^{(1)}(\tau) \pi^2(u^{1,o})(\tau,t) d\tau , \hst t \leq s \leq T. \label{pr2}
 \end{align}
We now introduce the following assumption regarding the measurability of admissible decisions ${\mathbb U}^{y^{i,u}}[0,T], i=1,2$.
\begin{assumptions}
\label{conj1}
Any admissible decentralized feedback information structure  $u^i \in {\mathbb U}^{y^{i,u}}[0,T]$ is adapted to $\Big\{{\cal G}_{0,t}^{I^{i,u}}: t \in [0,T]\Big\}$, $i=1,2$, and (\ref{ex50}) has a strong  $\Big\{{\cal G}_{0,t}^{I^{1,u}}: t \in [0,T]\Big\}-$
adapted solution $\pi^1(x)(\cdot)$ and (\ref{ex51}) has a strong  $\Big\{{\cal G}_{0,t}^{I^{2,u}}: t \in [0,T]\Big\}-$
adapted solution $\pi^2(x)(\cdot)$.
 \end{assumptions}

\noi Now, we can state the first main result.

\begin{theorem}
\label{dise}
Under the conditions of Assumptions~\ref{conj1} the optimal decisions $(u^{1,o}, u^{2,o})$ are given 
\begin{align}
u_t^{1,o} \equiv  -R_{11}^{-1}(t) B^{(1),*}(t) \Big\{K^1(t) \pi^1(x^o)(t) +r^1(t) \Big\}  -R_{11}^{-1}(t)R_{12}(t) \overline{u^{2,o}}(t),   \hst t \in [0,T],    \label{ex49i}
\end{align}
\begin{align}
u_t^{2,o}  \equiv  -R_{22}^{-1}(t) B^{(2),*}(t) \Big\{ K^2(t)    \pi^2(x^o)(t) +r^2(t)\Big\} -R_{22}^{-1}(t)R_{21}(t) \overline{u^{1,o}}(t), \hst t \in [0,T].  \label{ex49j}
\end{align}
where  $\pi^i(x^o)(\cdot), i=1,2$ satisfy the filter equations (\ref{ex50}), (\ref{ex51}), and  $\Big(K^i(\cdot), r^i(\cdot), \overline{x^o}(\cdot), \overline{u^{i,o}}(\cdot) \Big),i=1,2$ are solutions of  the ordinary differential equations (\ref{ex20pi}), (\ref{ex49f}), (\ref{ex49g}), (\ref{ex49h}), (\ref{ex101}), (\ref{ex106}).
\begin{align}
\dot{K}^i(t) &+ A^*(t) K^i(t)    + K^i(t) A(t)
 -K^i(t) B^{(i)}(t) R_{ii}^{-1}(t) B^{(i),*}(t) K^i(t)   \nonumber \\
 &+H(t)=0, \hst t \in [0,T), \hso i=1,2,   \label{ex20pi} \\
 K^i(T)&=M(T), \hso i=1,2 , \label{ex49f}
\end{align}
\begin{align}
\dot{r}^1(t)=&\Big\{  -A^{*}(t)  +\Phi^*(T,t) M(T)  \Psi_{K^1}(T,t) B^{(1)}(t) R_{11}^{-1}(t)B^{(1),*}(t) \nonumber \\                
&+ \Big(\int_{t}^T \Phi^*(s,t) H(s)  \Psi_{K^1}(s,t) ds \Big) B^{(1)}(t) R_{11}^{-1}(t) B^{(1),*}(t) \Big\} r^1(t)  \nonumber\\
&- \Big( \int_{t}^T  \Phi^*(s,t)  H(s)  \Psi_{K^1}(s, t) ds\Big) \Big(B^{(2)}(t) 
-B^{(1)}(t)R_{11}^{-1}(t) R_{12}(t) \Big)
\overline{u^{2,o}}(t), \nonumber \\
&-  \Phi^*(T,t)M(T) \Psi_{K^1}(T, t) \Big( B^{(2)}(t) -B^{(1)}(t)R_{11}^{-1}(t) R_{12}(t) \Big)  \overline{u^{2,o}}(t)   
   \hso t \in [0,T), \hso r^1(T)=0,   \label{ex49g} \\
\dot{r}^2(t)=& \Big\{  -A^{*}(t) +\Phi^*(T,t) M(T)  \Psi_{K^2}(T,t) B^{(2)}(t) R_{22}^{-1}(t)B^{(2),*}(t) \nonumber \\
&+ \Big(\int_{t}^T \Phi^*(s,t) H(s)  \Psi_{K^2}(s,t) ds \Big) B^{(2)}(t) R_{22}^{-1}(t) B^{(2),*}(t) \Big\}r^2(t) \nonumber \\
&- \Big( \int_{t}^T  \Phi^*(s,t)  H(s)  \Psi_{K^2}(s, t) ds\Big) \Big(B^{(1)}(t)  -B^{(2)}(t)R_{22}^{-1}(t) R_{21}(t) \Big)    \overline{u^{1,o}}(t)  \nonumber \\
&-  \Phi^*(T,t)M(T) \Psi_{K^2}(T, t) \Big(B^{(1)}(t)  -B^{(2)}(t)R_{22}^{-1}(t) R_{21}(t) \Big)    \overline{u^{1,o}}(t),
  \hso t \in [0,T), \hso r^2(T)=0,    \label{ex49h}
\end{align}  
\begin{align}
\dot{\overline{x^o}}(t)=A(t) \overline{x^o}(t) + B^{(1)}(t) \overline{u^{1,o}}(t)+ B^{(2)}(t) \overline{u^{2,o}}(t), \hst \overline{x^o}(0)=\overline{x}_0, \label{ex101}
\end{align}
\bea
   \left[ \begin{array}{c} \overline{u^{1,o}}(t) \\ \overline{u^{2,o}}(t) \end{array} \right] = -    \left[ \begin{array}{cc} I & R_{11}^{-1}(t)R_{12}(t)  \\ 
R_{22}^{-1}(t)R_{21}(t) & I \end{array} \right]^{-1}      \left[ \begin{array}{c} R_{11}^{-1}(t) B^{(1),*}(t) \Big\{K^1(t) \overline{x^o}(t) +r^1(t) \Big\} \\
R_{22}^{-1}(t) B^{(2),*}(t) \Big\{ K^2(t)    \overline{x^o}(t) +r^2(t)\Big\} \end{array} \right]. \label{ex106}
\eea

\end{theorem}
\begin{proof}  By invoking   Assumptions~\ref{conj1},  since $y^i(t) = \int_{0}^t  C^{[i]}(s)  \pi(x)(s)ds+ I^i(t)$, and $\pi^i(x)(\cdot)$ is a strong solution then ${\cal G}_{0,t}^{y^{i,u}} \subseteq {\cal G}_{0,t}^{I^{i,u}}$, $\forall t \in [0,T]$ and thus ${\cal G}_{0,t}^{y^{i,u}} = {\cal G}_{0,t}^{I^{i,u}}$ , $i=1, 2$. Hence, the optimality conditions of Theorem~\ref{theorem5.1o} are valid. 
Utilizing the independence of the innovations processes $I^1(\cdot)$ and $I^2(\cdot)$   then 
\begin{align}
\pi^1(u^2)(s,t) &={\mathbb E} \Big( u_s^2 | {\cal  G}_{0,t}^{y^1} \Big) ={\mathbb E} \Big( u_s^2\Big) \equiv \overline{u^2}(s), \hso t \leq s \leq T, \label{in2} \\ \pi^2(u^1)(s,t) &={\mathbb E} \Big( u_s^1 | {\cal  G}_{0,t}^{y^2} \Big) ={\mathbb E} \Big( u_s^1 \Big) \equiv \overline{u^1}(s), \hso t \leq s \leq T. \label{in2a}
\end{align}
Substituting (\ref{in2}), (\ref{in2a}) into (\ref{pr1}), (\ref{pr2}), and then  (\ref{pr1}), (\ref{pr2}) into (\ref{ex49}) we have 
\begin{align}
\pi^1(\psi^{o})(t)=& \Big\{ \Phi^*(T,t)M(T) \Psi_{K^1}(T,t) + \int_{t}^T \Phi^*(s,t)  H(s) \Psi_{K^1}(s,t)ds \Big\} \pi^1(x^o)(t) \nonumber \\
&+ \Phi^*(T,t)M(T) \int_{t}^T \Psi_{K^1}(T, \tau)\Big( B^{(2)}(\tau)-B^{(1)}(\tau)R_{11}^{-1}(\tau) R_{12}(\tau) \Big) \overline{u^{2,o}}(\tau) d\tau \nonumber \\ 
&+  \int_{t}^T \Phi^*(s,t)  H(s)  \int_{t}^s \Psi_{K^1}(s, \tau)\Big( B^{(2)}(\tau)  -B^{(1)}(\tau)R_{11}^{-1}(\tau) R_{12}(\tau) \Big) \overline{u^{2,o}}(\tau) d\tau ds \nonumber \\
&-\Phi^*(T,t) M(T) \int_{t}^T \Psi_{K^1}(T,\tau) B^{(1)}(\tau) R_{11}^{-1}(\tau)B^{(1),*}(\tau) r^1(\tau) d\tau \nonumber \\
&-\int_{t}^T \Phi^*(s,t) H(s) \int_{t}^s \Psi_{K^1}(s,\tau) B^{(1)}(\tau) R_{11}^{-1}(\tau) B^{(1),*}(\tau)r^1(\tau) d \tau ds,
 \label{ex49a}
\end{align}  
\begin{align}
\pi^2(\psi^{o})(t)=& \Big\{ \Phi^*(T,t)M(T) \Psi_{K^2}(T,t) + \int_{t}^T \Phi^*(s,t)  H(s) \Psi_{K^2}(s,t)ds \Big\} \pi^2(x^o)(t) \nonumber \\
&+ \Phi^*(T,t)M(T) \int_{t}^T \Psi_{K^2}(T, \tau) \Big( B^{(1)}(\tau)  
-B^{(2)}(\tau)R_{22}^{-1}(\tau) R_{21}(\tau) \Big) 
 \overline{u^{1,o}}(\tau) d\tau \nonumber \\ 
&+  \int_{t}^T \Phi^*(s,t)  H(s)  \int_{t}^s \Psi_{K^2}(s, \tau) \Big(B^{(1)}(\tau)   -B^{(2)}(\tau)R_{22}^{-1}(\tau) R_{21}(\tau) \Big)              \overline{u^{1,o}}(\tau) d\tau ds\nonumber \\
&-\Phi^*(T,t) M(T) \int_{t}^T \Psi_{K^2}(T,\tau) B^{(2)}(\tau) R_{22}^{-1}(\tau)B^{(2),*}(\tau) r^2(\tau) d\tau \nonumber \\
&-\int_{t}^T \Phi^*(s,t) H(s) \int_{t}^s \Psi_{K^2}(s,\tau) B^{(2)}(\tau) R_{22}^{-1}(\tau) B^{(2),*}(\tau)r^2(\tau) d \tau ds.
 \label{ex49b}
\end{align}  
Comparing (\ref{in54}) with the previous two equations then $K^i(\cdot), i=1,2$ are identified by the operators
\begin{align}
K^i(t)=  \Phi^*(T,t)M(T) \Psi_{K^i}(T,t) + \int_{t}^T \Phi^*(s,t)  H(s) \Psi_{K^i}(s,t)ds, \hst t \in [0,T], \hso i=1,2, \label{ex49c}
\end{align}
and $r^i(\cdot), i=1,2$ by the processes 
\begin{align}
r^1(t)=&\Phi^*(T,t)M(T) \int_{t}^T \Psi_{K^1}(T, \tau)\Big( B^{(2)}(\tau)-B^{(1)}(\tau)R_{11}^{-1}(\tau) R_{12}(\tau) \Big) \overline{u^{2,o}}(\tau) d\tau \nonumber \\ 
&+  \int_{t}^T \Phi^*(s,t)  H(s)  \int_{t}^s \Psi_{K^1}(s, \tau)\Big( B^{(2)}(\tau)  -B^{(1)}(\tau)R_{11}^{-1}(\tau) R_{12}(\tau) \Big) \overline{u^{2,o}}(\tau) d\tau ds \nonumber \\
&-\Phi^*(T,t) M(T) \int_{t}^T \Psi_{K^1}(T,\tau) B^{(1)}(\tau) R_{11}^{-1}(\tau)B^{(1),*}(\tau) r^1(\tau) d\tau \nonumber \\
&-\int_{t}^T \Phi^*(s,t) H(s) \int_{t}^s \Psi_{K^1}(s,\tau) B^{(1)}(\tau) R_{11}^{-1}(\tau) B^{(1),*}(\tau)r^1(\tau) d \tau ds,    \label{ex49d} 
\end{align}
\begin{align}
r^2(t)=& \Phi^*(T,t)M(T) \int_{t}^T \Psi_{K^2}(T, \tau) \Big( B^{(1)}(\tau)  
-B^{(2)}(\tau)R_{22}^{-1}(\tau) R_{21}(\tau) \Big) 
 \overline{u^{1,o}}(\tau) d\tau \nonumber \\ 
&+  \int_{t}^T \Phi^*(s,t)  H(s)  \int_{t}^s \Psi_{K^2}(s, \tau) \Big(B^{(1)}(\tau)   -B^{(2)}(\tau)R_{22}^{-1}(\tau) R_{21}(\tau) \Big)              \overline{u^{1,o}}(\tau) d\tau ds\nonumber \\
&-\Phi^*(T,t) M(T) \int_{t}^T \Psi_{K^2}(T,\tau) B^{(2)}(\tau) R_{22}^{-1}(\tau)B^{(2),*}(\tau) r^2(\tau) d\tau \nonumber \\
&-\int_{t}^T \Phi^*(s,t) H(s) \int_{t}^s \Psi_{K^2}(s,\tau) B^{(2)}(\tau) R_{22}^{-1}(\tau) B^{(2),*}(\tau)r^2(\tau) d \tau ds.  \label{ex49e}
\end{align}  
Differentiating both sides of (\ref{ex49c}) the operators $K^i(\cdot), i=1,2$ satisfy the following matrix differential equations (\ref{ex20pi}), (\ref{ex49f}). 
Differentiating both sides of (\ref{ex49d}), (\ref{ex49e})  the processes $r^i(\cdot), i=1,2$ satisfy the   differential equations (\ref{ex49g}), (\ref{ex49h}). Utilizing (\ref{in2}), (\ref{in2a}) we obtain the optimal strategies (\ref{ex49i}), (\ref{ex49j}).
Next, we determine $\overline{u^{i,o}}$ for $i=1,2$ from (\ref{ex49i}), (\ref{ex49j}). \\
Define the averages 
\begin{align}
\overline{x}(t) \tri {\mathbb E}\Big\{x(t)\Big\}= {\mathbb E} \Big\{ \pi^i(x)(t)  \Big\}, \hst i=1, 2.
\end{align}
Then $\overline{x^o}(\cdot)$ satisfies the ordinary differential equation (\ref{ex101}).
Taking the expectation of both sides of (\ref{ex49i}), (\ref{ex49j}) we deduce the corresponding  equations
\begin{align}
\overline{u^{1,o}}(t) &=  -R_{11}^{-1}(t) B^{(1),*}(t) \Big\{K^1(t) \overline{x^o}(t) +r^1(t) \Big\}  -R_{11}^{-1}(t)R_{12}(t) \overline{u^{2,o}}(t),   \hst t \in [0,T],    \label{ex102} \\
\overline{u^{2,o}}(t)  &=  -R_{22}^{-1}(t) B^{(2),*}(t) \Big\{ K^2(t)    \overline{x^o}(t) +r^2(t)\Big\} -R_{22}^{-1}(t)R_{21}(t) \overline{u^{1,o}}(t), \hst t \in [0,T].  \label{ex103}
\end{align}
The last two equations can be written in matrix form (\ref{ex106}). This completes the derivation.
\end{proof}

\noi Hence, the optimal strategies are computed from (\ref{ex49i}), (\ref{ex49j}), where the filter equations for $\pi^i(x^o)(\cdot), i=1,2$ satisfy (\ref{ex50}), (\ref{ex51}), while  $\Big(K^i(\cdot), r^i(\cdot), \overline{u^{i,o}}(\cdot), \overline{x^o}(\cdot) \Big),i=1,2$ are computed off-line utilizing  the ordinary differential equations (\ref{ex20pi}), (\ref{ex49f}), (\ref{ex49g}), (\ref{ex49h}), (\ref{ex101}), (\ref{ex106}). 

It is important to make the following observations.

\textbf{(O1):} The optimal strategies or laws (\ref{ex49i}), (\ref{ex49j}) are precisely the optimal strategies obtained in \cite{charalambous-ahmedFIS_Partii2012} for noiseless decentralized information structures. This property is analogous to that of optimal centralized strategies of fully and partially observed Linear-Quadratic-Gaussian systems. 

\textbf{(O2):} The filter equations for $\pi^i(x^o)(\cdot), i=1,2$ given by (\ref{ex50}), (\ref{ex51}) are nonlinear and may require higher order moments, leading to the so-called moment closure problem of nonlinear filtering. Furher analysis is required to determine whether the conditional error covariance in (\ref{ex50}), (\ref{ex51}) have the Kalman filter form. \\

Next, we state analogous results for distributed filtering problems.

\begin{corollary}
\label{disef}
Consider distributed filter dynamics (\ref{ex30}), (\ref{ex31}) with $B^{(i)}=0, i=1, 2$ and LQF pay-off (\ref{n3}). \\
Then the optimal strategies $(u^{1,o}, u^{2,o})$ are given by
\begin{align}
u_t^{i,o}=&-  R_{ii}^{-1}(t) \Big\{  m^i(t)+ \sum_{j=1}^2 E_{ij}(t) {\mathbb E} \Big( x^{j}(t) | {\cal G}_{0,t}^{y^{i}}\Big) + \sum_{j=1,j \neq i}^2  R_{ij}(t)  {\mathbb E} \Big( u_t^{j,o}  \Big)   \Big\}, \:  {\mathbb P}|_{ {\cal G}_{0,t}^{y^{i}}}-a.s., \: \forall i \in {\mathbb Z}_2, \label{n4ra}
\end{align}
where 
 $\widehat{x}^i(t)= Vector\Big\{ {\mathbb E} \{ x^1(t) | {\cal G}_{0,t}^{y^i}\}, {\mathbb E} \{ x^{2}(t) |{\cal G}_{0,t}^{y^i}\} \Big\},  i=1,2$ satisfy the linear Kalman filter equations
\begin{align}
d \widehat{x}^1(t) =& A(t) \widehat{x}^1(t)dt
+ P^1(t)   C^{[1],*}  D_{11}^{-1}(t)\Big(dy^{1}(t) - C^{[1]}(t) \widehat{x}^1(t) dt\Big), \hso \widehat{x}(0)=\bar{x}_0, \label{in50} \\
d\widehat{x}^2(t) =& A(t) \widehat{x}^2(t)dt +P^2(t)   C^{[2],*}  D_{22}^{-1}(t)\Big( dy^{2}(t)-C^{[2]}(t) \widehat{x}^2(t)dt\Big) , \hso \widehat{x^o}(0)
=\bar{x}_0, \label{in51} \\
\dot{P}^i(t) =& A(t)P^i(t)+ P^i(t)A(t) -P^2(t) C^{[i],*}(t) D_{ii}^{-1}(t) C^{[i]}(t)P^2(t)  \nonumber \\
&+G(t) G^{*}(t), \hso P^i(0)=P_0^i, \hso i=1,2. \label{in52}
\end{align}
and $\overline{u}(t) \tri Vector\Big\{ {\mathbb E} \{ u_t^1\}, {\mathbb E} \{ u_t^{2} \} \Big\}$ satisfy  the equations 
\begin{align}
 R(t)\overline{{ u}^o}(t) + E(t)\overline{x}(t) + m(t) =0,\hst  \frac{d}{dt} \overline{x}(t)=A \overline{x}(t), \hso \overline{x}(0)=\overline{x}_0. \label{fp1aa}
 \end{align}
\end{corollary}

\begin{proof} (\ref{n4ra}) is obtained from (\ref{n4r}) by setting $B^{(i)}=0, N=2$, and the discussion in  Section~\ref{static-feedback} (for filtering problems the observations and innovations generate the same filtrations). The filters (\ref{in50})-(\ref{in52}) are follow from the linear and Gaussian nature of the state and observation equations. Taking expectation of both sides of (\ref{n4ra}) yields (\ref{fp1aa}).
\end{proof}

Finally, we state a remark describing extensions of the previous examples.

\begin{remark}
\label{remPIP}
Theorem~\ref{dise} is easily generalized to the following arbitrary coupled  dynamics
\begin{align}
dx^i(t) =&A_{ii}(t)x^i(t) dt +  B^{(i)} u_t^i dt+ G_{ii} dW^i(t)  \nonumber \\
&+ \sum_{j=1, j\neq i}^N A_{ij}x^j(t) dt + \sum_{ j=1, j\neq i }^N B^{(j)}(t)u_t^j dt , \hso x^i(0)=x_0^i, \hso t \in (0,T], \hso i \in {\mathbb Z}_N \label{R1}
\end{align}
and information structures generated by observation equations with feedback
\bea
y^i = \int_{0}^t C_{ii}(s,y^i(s)) x(s) ds + \int_{0}^t D_{ii}^{\frac{1}{2}}(s)dB^i(s), \hso t \in [0,T], \hso i \in {\mathbb Z}_N. \label{R4}
\eea 
The optimal strategies are extensions of the ones given in Theorem~\ref{dise}. 
Similarly, one can generalize the filtering results of Corollary~\ref{disef} to the above  models with $B^{(i)}=0$.

\end{remark}

\section{Conclusions and Future Work}
\label{cf}
In this paper we have   considered team games for distributed stochastic differential decision systems, with decentralized noisy information patters for each DM, and we   derived  necessary and sufficient optimality conditions with respect to  team optimality and person-by-person optimality criteria, based on Stochastic Pontryagin's minimum principle.\\

 However, several additional issues remain to be investigated. Below, we provide a short list.

\begin{description}
\item[(F1)] In the derivation of optimality conditions we can relax some of the assumptions by considering spike or needle variations instead of strong variations of the decision strategies (or use relaxed strategies as in \cite{ahmed-charalambous2012a}). Moreover, for team games with non-convex action spaces ${\mathbb A}^i, i =1,2,\ldots, N$ and diffusion coefficients which depend on the decision variables it is necessary to derive optimality conditions based on second-order variations. 

\item[(F2)] The derivation of optimality conditions can be used in other type of games such as Nash-equilibrium games with decentralized noisy information structures for each DM, and minimax games.

\item[(F3)] It will be interesting determine whether  (\ref{ex50}) and (\ref{ex51}) are given by the Kalman filter equations.


\end{description}

\bibliographystyle{IEEEtran}
\bibliography{bibdata}

\begin{thebibliography}{10}
\providecommand{\url}[1]{#1}
\csname url@samestyle\endcsname
\providecommand{\newblock}{\relax}
\providecommand{\bibinfo}[2]{#2}
\providecommand{\BIBentrySTDinterwordspacing}{\spaceskip=0pt\relax}
\providecommand{\BIBentryALTinterwordstretchfactor}{4}
\providecommand{\BIBentryALTinterwordspacing}{\spaceskip=\fontdimen2\font plus
\BIBentryALTinterwordstretchfactor\fontdimen3\font minus
  \fontdimen4\font\relax}
\providecommand{\BIBforeignlanguage}[2]{{%
\expandafter\ifx\csname l@#1\endcsname\relax
\typeout{** WARNING: IEEEtran.bst: No hyphenation pattern has been}%
\typeout{** loaded for the language `#1'. Using the pattern for}%
\typeout{** the default language instead.}%
\else
\language=\csname l@#1\endcsname
\fi
#2}}
\providecommand{\BIBdecl}{\relax}
\BIBdecl

\bibitem{charalambous-ahmedFIS_Parti2012}
\BIBentryALTinterwordspacing
C.~D. Charalambous and N.~U. Ahmed, ``Centralized versus decentralized team
  games of distributed stochastic differential decision systems with noiseless
  information structures-{P}art {I}: General theory,'' \emph{Submitted to IEEE
  Transactions on Automatic Control}, p.~39, February 2013. [Online].
  Available: \url{http://arxiv.org/abs/1302.3452}
\BIBentrySTDinterwordspacing

\bibitem{ahmed-charalambous2012a}
\BIBentryALTinterwordspacing
N.~U. Ahmed and C.~D. Charalambous, ``Stochastic minimum principle for
  partially observed systems subject to continuous and jump diffusion processes
  and driven by relaxed controls,'' \emph{Submitted to SIAM Journal on Control
  and Optimization}, p.~23, June 2012. [Online]. Available:
  \url{http://arxiv.org/abs/1302.3455}
\BIBentrySTDinterwordspacing

\bibitem{fleming-rischel1975}
W.~Fleming and R.~Rischel, \emph{Deterministic and Stochastic Optimal
  Control}.\hskip 1em plus 0.5em minus 0.4em\relax Springer Verlag, 1975.

\bibitem{elliott1977}
R.~J. Elliott, ``The optimal control of stochastic system,'' \emph{SIAM Journal
  on Control and Optimization}, vol.~15, no.~5, pp. 756--778, 1977.

\bibitem{bismut1978}
J.~M. Bismut, ``An introductory approach to duality in optimal stochastic
  control,'' \emph{SIAM Review}, vol.~30, pp. 62--78, 1978.

\bibitem{elliott1982}
R.~J. Elliott, \emph{Stochastic Calculus and Applications}.\hskip 1em plus
  0.5em minus 0.4em\relax Springer-Verlag, 1982.

\bibitem{elliott-kohlmann1994}
R.~J. Elliott and M.~Kohlmann, ``The second order minimum principle and adjoint
  process,'' \emph{Stochastics \& Stochastic Reports}, vol.~46, pp. 25--39,
  1994.

\bibitem{peng1990}
S.~Peng, ``A general stochastic maximum principle for optimal control
  problems,,'' \emph{SIAM Journal on Control and Optimization}, vol.~28, no.~4,
  pp. 966--979, 1990.

\bibitem{yong-zhou1999}
J.~Yong and X.~Y. Zhou, \emph{Stochastic Controls, Hamiltonian Systems and
  {HJB} Equations}.\hskip 1em plus 0.5em minus 0.4em\relax Springer-Verlag,
  1999.

\bibitem{bensoussan1992a}
A.~Bensoussan, \emph{Stochastic Control of Partially Observable Systems}.\hskip
  1em plus 0.5em minus 0.4em\relax Cambridge University Press, 1982.

\bibitem{charalambous-hibey1996}
C.~D. Charalambous and J.~L. Hibey, ``Minimum principle for partially
  observable nonlinear risk-sensitive control problems using measure-valued
  decompositions,'' \emph{Stochastics \& Stochastic Reports}, pp. 247--288,
  1996.

\bibitem{ahmed-charalambous2007}
N.~U. Ahmed and C.~D. Charalambous, ``Minimax games for stochastic systems
  subject to relative entropy uncertainty: Applications to {SDE}'s on {H}ilbert
  {S}paces,'' \emph{Journal of Mathematics of Control, Signals and System},
  vol.~19, pp. 197--216, 2007.

\bibitem{varaiya-walrand1978}
P.~Varaiya and J.~Walrand, ``On delay sharing patterns,'' \emph{IEEE
  Transactions on Automatic Control}, vol.~23, no.~3, pp. 443--445, 1978.

\bibitem{witsenhausen1968}
H.~S. Witsenhausen, ``A counter example in stochastic optimum control,''
  \emph{SIAM Journal on Control and Optimization}, vol.~6, no.~1, pp. 131--147,
  1968.

\bibitem{witsenhausen1971}
------, ``Separation of estimation and control for discrete time systems,'' in
  \emph{Proceedinfs of the IEEE}, 1971, pp. 1557--1566.

\bibitem{ho-chu1972}
Y.-C. Ho and K.-C. Chu, ``Team decision theory and information structures in
  optimal control problems-part i,'' \emph{IEEE Transactions on Automatic
  Control}, vol.~17, no.~1, pp. 15--22, 1972.

\bibitem{kurtaran-sivan1973}
B.-Z. Kurtaran and R.~Sivan, ``Linear-{Q}uadratic-{G}aussian control with
  one-step-delay sharing pattern,'' \emph{IEEE Transactions on Automatic
  Control}, pp. 571--574, 1974.

\bibitem{sandell-athans1974}
N.~R. Sandell and M.~Athans, ``Solution of some nonclassical {LQG} stochastic
  decision problems,'' \emph{IEEE Transactions on Automatic Control}, vol.~19,
  no.~2, pp. 108--116, 1974.

\bibitem{kurtaran1975}
B.-Z. Kurtaran, ``A concice derivation of the {LQG} one-step-delay sharing
  problem solution,'' \emph{IEEE Transactions on Automatic Control}, vol.~20,
  no.~6, pp. 808--810, 1975.

\bibitem{ho1980}
Y.~Ho, ``Team decision theory and information structures,'' \emph{Proceedings
  of IEEE}, vol.~68, pp. 644--655, 1980.

\bibitem{bagghi-basar1980}
A.~Bagghi and T.~Basar, ``Teams decision theory for linear continuous-time
  systems,'' \emph{IEEE Transactions on Automatic Control}, vol.~25, no.~6, pp.
  1154--1161, 1980.

\bibitem{krainak-speyer-marcus1982a}
J.~Krainak, J.~L. Speyer, and S.~I. Marcus, ``Static team problems-part {I}:
  Sufficient conditions and the exponential cost criterion,'' \emph{IEEE
  Transactions on Automatic Control}, vol.~27, no.~4, pp. 839--848, 1982.

\bibitem{krainak-speyer-marcus1982b}
------, ``Static team problems-part {II}: Affine control laws, projections,
  algorithms, and the {LEGT} problem,'' \emph{IEEE Transactions on Automatic
  Control}, vol.~27, no.~4, pp. 848--859, 1982.

\bibitem{bansal-basar1987}
R.~Bansar and T.~Basar, ``Stochastic teams with nonclassical information
  revisited: When is an affine law optimal,'' \emph{IEEE Transactions on
  Automatic Control}, vol.~32, no.~6, pp. 554--559, 1987.

\bibitem{waal-vanschuppen2000}
P.~R. Wall and J.~H. van Schuppen, ``A class of team problems with discrete
  action spaces: Optimality conditions based on multimodularity,'' \emph{SIAM
  Journal on Control and Optimization}, vol.~38, no.~3, pp. 875--892, 2000.

\bibitem{bamieh-voulgaris2005}
B.~Bamieh and P.~Voulgaris, ``A convex characterization of distributed control
  problems in spatially invariant systems with communication constraints,''
  \emph{Systems and Control Letters}, vol.~54, no.~6, pp. 575--583, 2005.

\bibitem{aicardi-davoli-minciardi1987}
M.~Aicardi, F.~Davoli, and R.~Minciardi, ``Decentralized optimal control of
  markov chains with a common past information,'' \emph{IEEE Transactions on
  Automatic Control}, vol.~32, no.~11, pp. 1028--1031, 1987.

\bibitem{nayyar-mahajan-teneketzis2011}
A.~Nayyar, A.~Mahajan, and D.~Teneketzis, ``Optimal control strategies in
  delayed sharing information structures,'' \emph{IEEE Transactions on
  Automatic Control}, vol.~56, no.~7, pp. 1606--1620, 2011.

\bibitem{vanschuppen2011}
J.~H. van Schuppen, ``Control of distributed stochastic systems-introduction,
  problems, and approaches,'' in \emph{International Proceedings of the IFAC
  World Congress}, 2011.

\bibitem{lessard-lall2011}
L.~Lessard and S.~Lall, ``A state-space solution to the two-player optimal
  control problems,'' in \emph{Proceedings of 49th Annual Allerton Conference
  on Communication, Control and Computing}, 2011.

\bibitem{mahajan-martins-rotkowitz-yuksel2012}
A.~Mahajan, N.~Martins, M.~Rotkowitz, and S.~Yuksel, ``Information structures
  in optimal decentralized control,'' in \emph{In Proceedings of the 51st
  Conference on Decision and Control}, 2012.

\bibitem{gattami-bernhardsson-rantzer2012}
A.~Gattami, B.~M. Bernhardsson, and A.~Rantzer, ``Robust team decision
  theory,'' \emph{IEEE Transactions on Automatic Control}, vol.~57, no.~3, pp.
  794--798, 2012.

\bibitem{mishra-langbort-dullerud2012}
A.~Mishra, C.~Langbort, and G.~Dullerud, ``A team theoretic approach to
  decentralized control of systems with stochastic parameters,'' in \emph{In
  Proceedings of the 51st Conference on Decision and Control}, 2012, pp.
  2116--2121.

\bibitem{farokhi-johansson2012}
F.~Farokhi and K.~Johansson, ``Limited model information control design for
  linear discrete-time systems with stochastic parameters,'' in \emph{In
  Proceedings of the 51st Conference on Decision and Control}, 2012, pp.
  855--861.

\bibitem{marschak1955}
J.~Marschak, ``Elements for a theory of teams,'' \emph{Management Science},
  vol.~1, no.~2, 1955.

\bibitem{radner1962}
R.~Radner, ``Team decision problems,'' \emph{The Annals of Mathematical
  Statistics}, vol.~33, no.~3, pp. 857--881, 1962.

\bibitem{marschak-radner1972}
J.~Marschak and R.~Radner, \emph{Economic Theory of Teams}.\hskip 1em plus
  0.5em minus 0.4em\relax New Haven: Yale University Pres, 1972.

\bibitem{liptser-shiryayev1977}
R.~Liptser and A.~Shiryayev, \emph{Statistics of Random Processes Vol.1}.\hskip
  1em plus 0.5em minus 0.4em\relax Springer-Verlag New York, 1977.

\bibitem{charalambous-ahmedFIS_Partii2012}
\BIBentryALTinterwordspacing
C.~D. Charalambous and N.~U. Ahmed, ``Centralized versus decentralized team
  games of distributed stochastic differential decision systems with noiseless
  information structures-{P}art {II}: Applications,'' \emph{Submitted to IEEE
  Transactions on Automatic Control}, p.~39, February 2013. [Online].
  Available: \url{http://arxiv.org/abs/1302.3416}
\BIBentrySTDinterwordspacing

\end{thebibliography}

\end{document}